\documentclass[conference]{IEEEtran}

\usepackage{graphicx,bm}
\usepackage{stfloats}
\usepackage{amsmath}
\usepackage{amsfonts}
\usepackage{color}
\usepackage{amssymb}
\usepackage{amsthm}
\usepackage{caption2}
\usepackage{mathrsfs}
\usepackage{bm}
\usepackage{fancyhdr}
\usepackage{epstopdf}
\usepackage{txfonts}
\usepackage{dsfont}
\usepackage{multicol}

\newcounter{bla}
\begin{document}

\title{Generalized Cramer-Rao Bound for Joint Estimation
of Target Position and Velocity for {Active and} Passive  Radar
{Networks}}
\author{Qian He, \IEEEmembership{Member,~IEEE}, Jianbin Hu, Rick S. Blum,
\IEEEmembership{Fellow,~IEEE}, and Yonggang Wu
\thanks{The work of Q. He, J. Hu, and Yonggang Wu was supported by
the National Nature Science Foundation of China under Grants No.
61102142 and 61371184, the International Science and Technology
Cooperation and Exchange Research Program of Sichuan Province under
Grant No. 2013HH0006, and the Fundamental Research Funds for the
Central Universities under Grant No. ZYGX2013J015. The work of R. S.
Blum was supported by the National Science Foundation under Grant
No. ECCS-1405579.}
\thanks{Q. He, J. Hu, and Yonggang Wu are with University of Electronic
Science and Technology of China, Chengdu, Sichuan 611731 China
(emails: qianhe@uestc.edu.cn, 1205858265@qq.com, wyguestc@163.com).
R. S. Blum is with Lehigh University, Bethlehem, PA 18015 USA
(email: rblum@eecs.lehigh.edu). }} \IEEEoverridecommandlockouts

\maketitle

\begin{abstract}
In this paper, we derive the Cramer-Rao bound (CRB) for joint target
position and velocity estimation using an active or passive
distributed radar network under more general, and practically
occurring, conditions than assumed in previous work.  In particular,
the presented results allow nonorthogonal signals, spatially
dependent Gaussian reflection coefficients, and spatially dependent
Gaussian clutter-plus-noise.  These bounds allow designers to
compare the performance of their developed approaches, which are
deemed to be of acceptable complexity, to the best achievable
performance. If their developed approaches lead to performance close
to the bounds, these developed approaches can be deemed ``good
enough".  A particular recent study where algorithms have been
developed for a practical radar application which must involve
nonorthognal signals, for which the best performance is unknown,  is
a great example. The presented results in our paper do not make any
assumptions about the approximate location of the target being known
from previous target detection signal processing. In addition, for
{situations} in which we do not know some parameters accurately, we
also derive the mismatched CRB. Numerical investigations of the mean
squared error of the maximum likelihood estimation are employed to
support the validity of the CRBs. In order to demonstrate the
utility of the provided results to a topic of great current
interest, the numerical results focus on a passive radar system
using the Global System for Mobile communication (GSM) cellar
system.
\end{abstract}




\IEEEpeerreviewmaketitle

\section{Introduction}

The focus of this paper is on new Cramer-Rao bounds (CRB) for
estimation of target position and velocity from distributed radar
networks, sometimes called MIMO radar systems or multistatic radar
systems
\cite{Visa14,PuWang13,YaoYu,Alex,XiufengSong,Yimin12,Hana,JianLi07},
operating under more general, and practically occurring, conditions
than assumed in previous work. In particular, the presented results
allow nonorthogonal signals, spatially dependent Gaussian reflection
coefficients, and spatially dependent Gaussian clutter-plus-noise,
which are cases of great practical interest.  In fact, one could
argue that all of these conditions are true in any real system. The
initial CRBs we present are applicable to both active and passive
radar systems, provided the signals of opportunity in the passive
systems are assumed to be perfectly estimated from, for example, the
direct path reception.  These results further assume all the
parameters of the observations model are known, including the
covariance matrices of the zero-mean Gaussian noise-plus-clutter and
reflection coefficients.  The mismatched CRB results given in this
paper even allow cases where the model assumed by the estimation
algorithm is incorrect, including cases where the model for the
direct path signal may involve unmodeled noise, interference, or
some other imperfection. Similarly, the reflection coefficients,
noise, andor interference may be incorrectly modeled and the
mistmatched CRB will still provide a lower bound on perfromance.

While both passive and active radar systems are of great interest,
passive radar systems may have attracted even greater attention over
the past few years due to the tremendous advantages they provide
from using existing communication signals to implement a radar,
essentially borrowing the already existing transmitter
infrastructure and providing no electronic evidence that a radar is
operating in a given area. Passive radar, as the name implies, is a
radar system which receives only. Instead of actively transmitting
signals, it works passively by gathering signals from
non-cooperative illuminators of opportunity and reflections from
objects in the monitored area to make decisions or provide
information about targets. Since transmitters are not required in a
passive radar, it has the advantages of low implementation costs,
stealth, and the ability to operate in a wide frequency band without
concerns of causing interference to existing wireless systems. {For
these reasons, passive radar systems have attracted the attention
of} the international radar community. Passive radar {systems} based
on FM \cite{Howland:2005} and digital illuminators ({DAB/DVB-T})
\cite{Glende:2007}, or satellite-borne illuminators
\cite{Griffiths:2002}, WIFI \cite{Guo:2008} and {Global System for
Mobile communication} (GSM) \cite{Tan:2005} signals have been
{previously investigated mainly from prototypes or measurements or
very simple analytical models}. The factors that affect the
detection performance of passive coherent location radar systems are
{discussed} in \cite{Griffiths:2005}. The ambiguity functions of a
set of off-air measurements of signals that may be used for passive
coherent location (PCL) radar systems are presented and analyzed in
\cite{Griffiths and Baker:2005 }. The problem of target detection in
passive MIMO radar (PMR) networks comprised of non-cooperative
transmitters and multichannels is addressed in \cite{Hack:2014}.


As described later, the CRB is a lower bound, in a certain sense, on
the covariance matrix of all unbiased estimators.  It is a useful
tool for evaluating the best possible estimation performance of a
radar system. {A derivation of the stochastic CRB is provided in
\cite{Stoica:2001}.} The CRB expressions for the estimation of range (time delay),
velocity (Doppler shift) and direction of a point target using an
active radar or sonar array are given in \cite{Dogandzic:2001}. The
CRB of DOA estimation of a non-stationary target for a MIMO radar
with colocated antennas for a general time division multiplexing
(TDM) scheme is computed in \cite{Rambach:2013}. The CRB for
bistatic radar channels {is derived in \cite{Greco:2011}, which also
exploits the relationship} between the ambiguity function and the
CRB. {Cramer-Rao-like bounds for the estimation of a deterministic
parameter in the presence of random nuisance parameters are derived
in \cite{Gini:2012}, \cite{RW:1978}.}

{For the case of} multiple transmit and receive antennas {employed}
in a distributed active radar setting, {\cite{He:2010} describes
the} CRB under the assumption of orthogonal signals, spatially
independent reflection coefficients, and spatially independent
clutter-plus-noise. For estimation of the position and velocity of a
single target using a passive radar, the CRB and ambiguity functions
are considered in \cite{Stinco:2012} for a multiple transmitter and
receiver radar, but only for the case where a single transmitter and
receiver pair is selected from among a much larger set of possible
pairs.  This work does not consider the effect of signal
nonorthogonality or spatially dependent reflection coefficients or
noise. Under the same assumptions employed in \cite{He:2010}, the
CRB has been derived for passive radar settings with well estimated
signals of opportunity in \cite{He:2014,Gogineni:2014_2}.

Thus, none
of the published work has given the CRB for the important and
practical case of nonorthogonal signals, spatially dependent
reflection coefficients, and spatially dependent noise for joint
target position and velocity estimation performance using a
distributed passive or active radar network employing all signals
available from the multiple transmit and receive antenna paths
in an optimum manner.   This result is extremely useful
since it describes the best achievable performance for some
important cases for the first time.  Knowing this best achievable
performance allows designers to compare the performance of their
developed approaches, to these bounds.   If their developed
approaches lead to performance close to the bounds, these developed
approaches can be deemed "good enough" while these developed
approaches are typically constrained to have acceptable complexity.
The very recent work in \cite{addnonorth} provides an excellent
example where these results can be extremely useful. In
\cite{addnonorth}, a very practical scenario is considered where a
number of transmitters of opportunity send digital TV signals that
can not be accurately modeled by assuming the transmitters send a
set of nonorthogonal signals.  The work in \cite{addnonorth} presents
an interesting suboptimum algorithm for implementing a radar
employing these nonorthogonal signals.  However, it is not known how
far the performance of the suggested approach in \cite{addnonorth}
is from the optimum achievable performance.  Such information would
be extremely useful in judging if the approach suggested in
\cite{addnonorth} provides a good tradeoff in terms of performance
and complexity.   Similar questions arise in many related practical
applications, some of which involve active radars.


In this paper, we consider these more general cases and derive a
generalized CRB and mismatched CRB for joint location and velocity
estimation in passive and active distributed radar networks. The
presented results do not assume the approximate location of the
target is known from previous target detection signal processing,
unlike the previous results employing optimum processing using all
available antennas \cite{He:2010,He:2014,Gogineni:2014_2}. A
closed-form Fisher information matrix {(FIM)} is presented. In a few
representative cases, the generalized or mismatched CRB is
numerically compared with the mean-squared error (MSE) from maximum
likelihood (ML) estimation to show consistency at higher
{signal-to-noise ratio} (SNR). We use GSM signals as illuminators
for our numerical passive radar investigations. The rest of this
paper is organized as follows. The signal model for active and
passive distributed radar networks is presented in Section
\ref{sec:SigMod}. The ML estimate is analyzed in Section
\ref{sec:MLE}. In Section \ref{sec:GCRB}, the generalized CRB is
derived. In Section \ref{mm}, we derive the mismathed CRB.
Performance analysis and numerical examples are presented in Section
\ref{sec:Result}. Finally, conclusions are drawn in Section
\ref{sec:Conlusion}.

Throughout this paper, the notation for transpose is $(.)^{\dag}$,
while that for complex conjugate is $(.)^H$. The symbol
$Diag\{\cdot\}$ denotes a block diagonal matrix with the matrices in
the braces being the diagonal blocks,
$\mathcal{CN}({{\bm{\mu}},{\bm{R}}})$ denotes a complex Gaussian
distribution with mean vector $ \bm{\mu} $ and covariance matrix $
{\bm{R}} $, ${\mathbb{E}_{{\bm{r}}|{\bm{\theta
}},{\bm{\alpha}}}}\left\{ \cdot \right\}$ implies taking expectation
with respect to the probability density function (pdf) $p\left(
{{\bm{r}}|{\bm{\theta }},{\bm{\alpha}}} \right)$, $\text{Tr}\left(
\cdot \right)$ denotes the trace of a matrix, $\otimes$ represents
the Kronecker product, $\Re(\cdot)$ means taking the real part,
$\odot$ represents the Hadamard product, and $\text{vec}(\cdot)$
denotes the column vectorizing operator which stacks the columns of
a matrix in a column vector.

\section{Signal Model}\label{sec:SigMod}

Consider a distributed radar system with $M$ widely spaced single
antenna transmit stations and $N$ widely spaced single antenna
receive stations, located at $(x_m^t, y_m^t)$, $m=1,\ldots,M$ and $(
x_n^r, y_n^r)$, $n=1,\ldots,N$ in a two-dimensional Cartesian
coordinate system, respectively. The lowpass equivalent time-sampled
version of the signal transmitted from the $m$th transmit station at
time instant $k{T_s}$ is $\sqrt {{E_m}} {s_m}\left( k,\bm{\alpha_m }
\right)$, where $T_s$ is the sampling period, $k$ ($k=1,\ldots,K$)
is an index running over the different time samples, $\bm{\alpha_m
}$ denotes a vector of parameters needed to describe the waveform,
and the waveform is normalized using $\sum\nolimits_{k = 1}^{K}
|{{s_m}\left( k,\bm{\alpha_m }\right)}|^2 T_s= 1$.  Let ${E_m}$
denote {the} energy transmitted by the ${m}$th transmit antenna.
Then the received waveform at the ${n}$th receiver at time $kT_s$ is
\begin{equation}\label{rnm}
{r_n}\left( k \right) = \sum\limits_{m = 1}^M {\sqrt{\frac{
{{E_m}{P_0}} }{{d^2_{tm}d^2_{rn}}}}{\zeta _{nm}}{s_m}\left( {k{T_s}
- {\tau _{nm}},{\bm{\alpha} _m}} \right){e^{j2\pi {f_{nm}}k{T_s}}} +
{w_n}\left( k \right)},
\end{equation}
where ${\tau _{nm}}$, $f_{nm}$, and ${\zeta_{nm}}$ represent the
time delay, Doppler shift, and reflection coefficient corresponding
to the $nm$th path, respectively. {The variable} $d_{tm}$ denotes
the distance between the target and the ${m}$th transmitter, {while}
$d_{rn}$ denotes the distance between the target and the ${n}$th
receiver. {The term ${w_n}\left( k \right)$ denotes
clutter-plus-noise at the ${n}$th receiver at {time} $kT_s$.} The
received signal strength at $d_{tm}$=$d_{rn}$=1 is $\sqrt{E_mP_0}$,
so $P_0$ {denotes the} ratio of received energy at
$d_{tm}$=$d_{rn}$=1 to transmitted energy. The reflection
coefficient ${\zeta_{nm}}$ is assumed to be constant over the
observation interval {and to have} a known complex Gaussian
statistical model \cite{handbook}. Assume the position $\left( {x,y}
\right)$ and velocity $(v_x, v_y)$ of the target are deterministic
{unknowns}. The distances $d_{tm}$ and $d_{rn}$ {are expressed in
terms of $\left( {x,y} \right)$ as}
\begin{align}
d_{tm} = \sqrt {{{\left( {x_m^t - x} \right)}^2} + {{\left( {y_m^t -
y} \right)}^2}},
\end{align}
\begin{align}
d_{rn} = \sqrt {{{\left( {x_n^r - x} \right)}^2} + {{\left( {y_n^r -
y} \right)}^2}}.
\end{align}
The time delay ${\tau _{nm}}$ is {also} a function of the unknown
target position $\left( {x,y} \right)$
\begin{align}
  {\tau _{nm}} &= \frac{{\sqrt {{{\left( {x_m^t - x} \right)}^2}
  + {{\left( {y_m^t - y} \right)}^2}}  + \sqrt {{{\left( {x_n^r
  - x} \right)}^2} + {{\left( {y_n^r - y} \right)}^2}} }}{c}\notag\\
   &= \frac{{d_{tm} + d_{rn}}}{c}\label{taonm},
\end{align}
where $c$ denotes the speed of light, The Doppler shift ${f_{nm}}$
is a function of the unknown target position $(x, y)$ and velocity
$(v_x, v_y)$ {given by}
\begin{align}
{f_{nm}} = \frac{{{v_x}\left( {x_m^t - x} \right) + {v_y}\left(
{y_m^t - y} \right)}}{{\lambda d_{tm}}}+ \frac{{{v_x}\left( {x_n^r -
x} \right) + {v_y}\left( {y_n^r - y} \right)}}{{\lambda
d_{rn}}}\label{fnm},
\end{align}
where $\lambda$ denotes the wavelength. Define an unknown parameter
vector ${\bm{\theta }}$ that collects the parameters to be estimated
\begin{equation}
{\bm{\theta }} = {\left[ {x,y,{v_x},{v_y}} \right]^\dag}.
\end{equation}
The observations from the $n$th receiver can be {expressed} as
\begin{align}
{\bm{r}}_n &= \left[ {{{r}_n}\left( 1 \right),{{r}_n}\left( 2 \right), \cdots ,{{r}_n}\left( K \right)} \right]^\dag \\
&= {{\bm{U}}}_n {{\bm{\zeta }_n}} + {{\bm{w}}}_ n,
\end{align}
where ${\bm{U}}_n$ is a $K \times M$ matrix that collects the time
delayed and Doppler shifted signals at the $n$th receiver {as}
\begin{align}
 {\bm{U}}_n^{} &= \left[{{\bm{u}}_n\left( 1 \right)},{{\bm{u}}_n\left( 2 \right)},
   \dots,{{\bm{u}}_n\left( K \right)}\right]^\dag,
\end{align}
where
\begin{equation}
{{\bm{u}}_n}\left( k \right) = {\left[ {{u_{n1}}\left( k
\right),{u_{n2}}\left( k \right), \cdots ,{u_{nM}}\left( k \right)}
\right]^\dag},
\end{equation}
and
\begin{equation}
{u_{nm}}(k) = {\sqrt{\frac{ {{E_m}{P_0}}
}{{d^2_{tm}d^2_{rn}}}}}{s_m}(k{T_s} - {\tau _{nm}},{\bm{\alpha
_m}}){e^{j2\pi {f_{nm}}k{T_s}}}.
\end{equation}
The $M\times 1$ reflection coefficient vector ${\bm{\zeta }_n}$ can
be expressed as ${\bm{\zeta }_n} = {\left[ {{\zeta _{n1}}, \cdots
,{\zeta _{nM}}} \right]^\dag}$. Denote the vector of noise samples
at the $n$th receiver as ${\bm{w}}_n= {\left[ {{w_n}\left( 1
\right), \cdots ,{w_n}\left( K \right)} \right]^\dag}$. The
observations from {the set of all receivers} can be written as
\begin{align}
  {{\bm{r}}}=& \left[ {{\bm{r}}^\dag_1,\bm{r}^\dag_2, \cdots ,\bm{r}^\dag_N } \right]^\dag\notag\\
=& {\bm{S\bm\zeta }} + {\bm{w}}\label{rsignal},
\end{align}
where $\bm{S}$ collects the time delayed and Doppler shifted signals
from all paths
\begin{align}
{\bm{S}} = Diag\{ {{\bm{U}}_1},{{\bm{U}}_2}, \ldots ,{{\bm{U}}_N}\}.
\end{align}
The $\bm{\zeta }$ in (\ref{rsignal}) collects reflection
coefficients for all paths
\begin{equation}
{\bm{\zeta }} = {\left[ {{\bm{\zeta} _{1}^\dag}, \cdots
,{\bm{\zeta}_{N}^\dag}} \right]^\dag},
\end{equation}
and it is assumed that $\bm{\zeta }$ is a complex Gaussian
distributed {vector} with {zero mean} and covariance matrix
${{\bm{R}}} = \mathbb{E}\{ \bm{\zeta}\bm{\zeta}^H \}$, i.e.
${\bm{\zeta }} \sim\mathcal{CN}\left( {{\bm{0}},{\bm{R}}} \right)$.
The $\bm{w}$ in (\ref{rsignal}) denotes the clutter-plus-noise
vector
\begin{equation}
{\bm{w}} = {\left[ {\bm{w}^\dag_1, \cdots ,{\bm{w}}^\dag_N }
\right]^\dag},
\end{equation}
which is assumed to be complex Gaussian distributed with zero mean
and covariance matrix ${\bm{Q}}=\mathbb{E}\{\bm{w}\bm{w}^H\}$, i.e.,
${\bm{w}} \sim \mathcal{CN}\left( {{\bm{0}},{\bm{Q}}} \right)$.
Assume that the noise vector $\bm w$ is independent from the
reflection coefficient vector $\bm \zeta$.

\subsection{{Maximum Likelihood Estimation}}\label{sec:MLE}

In this and the next section (Sections II and III), we assume
$\bm{S}$ (and thus ${\bm{\alpha }}$), $\bm{Q}$, and $\bm{R}$ are
known to the estimation algorithm. We address other cases later.
Using the signal model in (\ref{rsignal}) and the fact that the
linear combination of two Gaussian vectors {is} also Gaussian, the
likelihood function conditioned on the waveform parameter vector
\begin{equation}
{\bm{\alpha }} = {\left[ {\bm{\alpha_1 },\ldots,\bm{\alpha_M }}
\right]^\dag},
\end{equation}
can be obtained as
\begin{equation}
p\left( {{\bm{r}}|{\bm{\theta }},\bm{\alpha}} \right) =
\frac{1}{{{\pi ^{KN}}\det ( {\bm{C}} )}}\exp ( { -
{{\bm{r}}^H}{{\bm{C}}^{ - 1}}{\bm{r}}} )\label{Prtheta},
\end{equation}
where ${\bm{C}}$ denotes the covariance matrix
\begin{align}
{\bm{C}} &= \mathbb {E}\left\{ {\left( {{\bm{S\zeta }} + {\bm{w}}}
\right){{\left( {{\bm{S\zeta }} + {\bm{w}}} \right)}^H}} \right\} \notag\\
&= \mathbb {E}\left\{ {{\bm{S\zeta }}{{\bm{\zeta }}^H}{{\bm{S}}^H} +
{\bm{w}}{{\bm{w}}^H}} \right\} \notag\\ &= {\bm{SR}}{{\bm{S}}^H} +
{\bm{Q}}\label{C}.
\end{align}
 The log-likelihood function
can be written as
\begin{align}
  L\left( {{\bm{r}}|{\bm{\theta }},\bm{\alpha}} \right)
  &= \ln p\left( {{\bm{r}}|{\bm{\theta }},\bm{\alpha}} \right)\notag \\
   &=  - {{\bm{r}}^H}{{\bm{C}}^{ - 1}}{\bm{r}} - \ln \left( {\det
   \left( {\bm{C}} \right)} \right) - KN\ln \left( \pi  \right)\label{logFun}.
\end{align}

Neglecting the last constant term of the second line in
(\ref{logFun}) and assuming known or perfectly estimated
$\bm{\alpha}$, the (ML) estimate of the unknown parameter vector
$\bm{\theta}$ can be calculated as
\begin{align}
{{\bm{\hat \theta }}_{ML}} &=\arg \mathop {\max }\limits_{\bm{\theta}}{L(\bm{r}|\bm{\theta},\bm{\alpha})}\notag \\
&=\arg \mathop {\max }\limits_{\bm{\theta}} \left\{ { -
{{\bm{r}}^H}{{\bm{C}}^{ - 1}}{\bm{r}} - \ln \left( {\det \left(
{\bm{C}} \right)} \right)} \right\}.
\end{align}

\section{Generalized Cramer-Rao Bound }\label{sec:GCRB}

In this section, we provide the CRB for jointly estimating the
target location $\left( {x,y} \right)$ and velocity $\left(
{{v_x},{v_y}} \right)$ for the case where $\bm{S}$ (and thus
${\bm{\alpha }}$), $\bm{Q}$, and $\bm{R}$ are known to the
estimation algorithm. The first step in obtaining the CRB is to
compute the FIM, which is a $4\times4$ matrix related to the second
order derivatives of the log-likelihood function
 \begin{equation}
 \begin{split}
  {\bm{J}}\left( {\bm{\theta }|\bm{\alpha }} \right) &=
  {\mathbb{E}_{{\bm{r}}|{\bm{\theta }},{\bm{\alpha }}}}\left\{
  {{\nabla _{\bm{\theta }}}L\left( {{\bm{r}}|{\bm{\theta }},{\bm{\alpha }}} \right)
  {{\left[ {{\nabla _{\bm{\theta }}}L\left( {{\bm{r}}|{\bm{\theta }},{\bm{\alpha }}}
  \right)} \right]}^\dag}} \right\}.
 \end{split}
\end{equation}
Considering the likelihood is a function of ${\tau _{nm}}$,
${f_{nm}}$, ${d_{tm}}$, and ${d_{rn}}$ ($n=1, \cdots, N$,
$m=1,\cdots, M$), which depend on ${\bm{\theta }} = {\left[
{x,y,{v_x},{v_y}} \right]^\dag}$, we define an intermediate
parameter vector

\begin{align}
\bm{\vartheta}  = &{[{{\bm{\tau
}}^\dag},{{\bm{f}}^\dag},{\bm{d}}_t^\dag,{\bm{d}}_r^\dag]^\dag}\\\notag
=&[{\tau _{11}},{\tau _{12}}, \cdots ,{\tau
_{NM}},{f_{11}},{f_{12}}, \cdots ,{f_{NM}},\\\notag
&{d_{t1}},{d_{t2}}, \cdots ,{d_{tM}},{d_{r1}},{d_{r2}}, \cdots
{d_{rN}}{]^\dag}
\end{align}
where $\bm{\tau}=\left[ {\tau _{11}},{\tau _{12}}, \cdots ,{\tau
_{NM}}\right]^\dag$, $\bm{f}=\left[{f_{11}},{f_{12}}, \cdots
,{f_{NM}} \right]^\dag$, $\bm{d_{t}}=\left[ {d_{t1}},{d_{t2}},
\cdots ,{d_{tM}}\right]^\dag$ and  $\bm{d_{r}}=\left[
{d_{r1}},{d_{r2}}, \cdots ,{d_{rN}}\right]^\dag$ collect the unknown
time {delays}, {Doppler shifts}, and distance parameters,
respectively. According to the chain rule, the FIM can be derived by
\begin{equation}\label{Fisher}
{\bm{J}}\left( {\bm{\theta }|\bm{\alpha }} \right) = \left( {{\nabla
_{\bm{\theta }}}{\bm{\vartheta} ^\dag}} \right){\bm{J}}\left(
\bm{\vartheta}|\bm{\alpha }  \right){\left( {{\nabla _{\bm{\theta
}}}{\bm{\vartheta} ^\dag}} \right)^\dag},
\end{equation}
{where ${\bm{J}}\left( {{\bm{\vartheta}} |{\bm{\alpha }}} \right) =
{\mathbb{E}_{{\bm{r}}|{\bm{\vartheta}}, {\bm{\alpha }} }}\left\{ {{\nabla
_{\bm{\vartheta}} }L\left( {{\bm{r}}|{\bm{\vartheta}} ,{\bm{\alpha }}}
\right){{\left[ {{\nabla _{\bm{\vartheta}} }L\left( {{\bm{r}}|{\bm{\vartheta}}
,{\bm{\alpha }}} \right)} \right]}^\dag}} \right\}$.}

{\subsection{Calculation of ${\nabla _{\bm{\theta }}}{{\bm{\vartheta}}
^\dag}$}}

Recalling (6) and (22), we have
\begin{equation}
{\nabla _{\bm{\theta }}}{{\bm{\vartheta}} ^\dag} = \left[ \begin{array}{l}
{\bm{F}}{\kern 1pt} {\kern 1pt} {\kern 1pt} {\kern 1pt} {\kern 1pt}
{\kern 1pt} {\kern 1pt} {\kern 1pt} {\kern 1pt} {\bm{G}}{\kern 1pt}
{\kern 1pt} {\kern 1pt} {\kern 1pt} {\kern 1pt} {\kern 1pt} {\kern
1pt} {\kern 1pt} {\kern 1pt} {\kern 1pt} {{\bm{D}}_t}{\kern 1pt}
{\kern 1pt} {\kern 1pt} {\kern 1pt} {\kern 1pt} {\kern 1pt} {\kern
1pt} {\kern 1pt}
{\kern 1pt} {\kern 1pt} {{\bm{D}}_r}\\
{\bm{0}}{\kern 1pt} {\kern 1pt} {\kern 1pt} {\kern 1pt} {\kern 1pt}
{\kern 1pt} {\kern 1pt} {\kern 1pt} {\kern 1pt} {\kern 1pt}
{\bm{H}}{\kern 1pt} {\kern 1pt} {\kern 1pt} {\kern 1pt} {\kern 1pt}
{\kern 1pt} {\kern 1pt} {\kern 1pt} {\kern 1pt} {\kern 1pt} {\kern
1pt} {\bm{0}}{\kern 1pt} {\kern 1pt} {\kern 1pt} {\kern 1pt} {\kern
1pt} {\kern 1pt} {\kern 1pt} {\kern 1pt} {\kern 1pt} {\kern 1pt}
{\kern 1pt} {\kern 1pt} {\kern 1pt} {\kern 1pt} {\kern 1pt} {\kern
1pt} {\kern 1pt} {\kern 1pt} {\bm{0}}
\end{array} \right]\label{fgdt},
\end{equation}
where
\begin{equation}
{\bm{F}} = \left[ {\begin{array}{*{20}{c}}
  {\frac{{\partial {\tau _{11}}}}{{\partial x}}}&{\frac{{\partial {\tau _{12}}}}{{\partial x}}}& \cdots &{\frac{{\partial {\tau _{NM}}}}{{\partial x}}} \\
  {\frac{{\partial {\tau _{11}}}}{{\partial y}}}&{\frac{{\partial {\tau _{12}}}}{{\partial y}}}& \cdots &{\frac{{\partial {\tau _{NM}}}}{{\partial y}}}
\end{array}} \right],
\end{equation}
\begin{equation}
{\bm{G}} = \left[ {\begin{array}{*{20}{c}}
  {\frac{{\partial {f_{11}}}}{{\partial x}}}&{\frac{{\partial {f_{12}}}}{{\partial x}}}& \cdots &{\frac{{\partial {f_{NM}}}}{{\partial x}}} \\
  {\frac{{\partial {f_{11}}}}{{\partial y}}}&{\frac{{\partial {f_{12}}}}{{\partial y}}}& \cdots &{\frac{{\partial {f_{NM}}}}{{\partial y}}}
\end{array}} \right],
\end{equation}

\begin{equation}
{\bm{H}} = \left[ {\begin{array}{*{20}{c}}
  {\frac{{\partial {f_{11}}}}{{\partial {v_x}}}}&{\frac{{\partial {f_{12}}}}{{\partial {v_x}}}}& \cdots &{\frac{{\partial {f_{NM}}}}{{\partial {v_x}}}} \\
  {\frac{{\partial {f_{11}}}}{{\partial {v_y}}}}&{\frac{{\partial {f_{12}}}}{{\partial {v_y}}}}& \cdots &{\frac{{\partial {f_{NM}}}}{{\partial {v_y}}}}
\end{array}} \right],
\end{equation}
\begin{equation}
{{\bm{D}}_t} = \left[ {\begin{array}{*{20}{c}}
{\frac{{\partial {d_{t1}}}}{{\partial x}}}&{\frac{{\partial {d_{t2}}}}{{\partial x}}}& \cdots &{\frac{{\partial {d_{tM}}}}{{\partial x}}}\\
{\frac{{\partial {d_{t1}}}}{{\partial y}}}&{\frac{{\partial
{d_{t2}}}}{{\partial y}}}& \cdots &{\frac{{\partial
{d_{tM}}}}{{\partial y}}},
\end{array}{\kern 1pt} {\kern 1pt} {\kern 1pt} {\kern 1pt} {\kern 1pt} {\kern 1pt} {\kern 1pt} {\kern 1pt} {\kern 1pt} {\kern 1pt}
} \right],
\end{equation}
and
\begin{equation}
{{\bm{D}}_r} = {\kern 1pt} \left[ {{\kern 1pt} {\kern 1pt} {\kern
1pt} \begin{array}{*{20}{c}}
{\frac{{\partial {d_{r1}}}}{{\partial x}}}&{\frac{{\partial {d_{r2}}}}{{\partial x}}}& \cdots &{\frac{{\partial {d_{rN}}}}{{\partial x}}}\\
{\frac{{\partial {d_{r1}}}}{{\partial y}}}&{\frac{{\partial
{d_{r2}}}}{{\partial y}}}& \cdots &{\frac{{\partial
{d_{rN}}}}{{\partial y}}}
\end{array}} \right].
\end{equation}
Using calculations drawing on (2)-(\ref{fnm}), the elements of the
matrices in (25)-(29) will be described as
\begin{equation}
{a_{nm}} = \frac{{\partial {\tau _{nm}}}}{{\partial x}} =
\frac{1}{c}\left( {\frac{{x - x_m^t}}{{d_{tm}}} + \frac{{x -
x_n^r}}{{d_{rn}}}} \right),
\end{equation}
\begin{equation}\label{}
{b_{nm}} = \frac{{\partial {\tau _{nm}}}}{{\partial y}} =
\frac{1}{c}\left( {\frac{{y - y_m^t}}{{d_{tm}}} + \frac{{y -
y_n^r}}{{d_{rn}}}} \right),
\end{equation}
\begin{align}\label{}
{e_{nm}} =& \frac{{\partial {f_{nm}}}}{{\partial x}} =  - \frac{{{v_x}}}{\lambda }\left( {\frac{1}{{d_{tm}}} + \frac{1}{{d_{rn}}}} \right)\notag \\
&+ \frac{{\left( {x_m^t - x} \right)}}{{\lambda {{\left( {d_{tm}} \right)}^3}}}\left[ {{v_x}\left( {x_m^t - x} \right) + {v_y}\left( {y_m^t - y} \right)} \right]\notag\\
&+ \frac{{\left( {x_n^r - x} \right)}}{{\lambda {{\left( {d_{rn}}
\right)}^3}}}\left[ {{v_x}\left( {x_n^r - x} \right) + {v_y}\left(
{y_n^r - y} \right)} \right],
\end{align}
\begin{align}\label{}
 {g_{nm}} &= \frac{{\partial {f_{nm}}}}{{\partial y}} =  - \frac{{{v_y}}}{\lambda }\left( {\frac{1}{{d_{tm}}} + \frac{1}{{d_{rn}}}} \right)\notag \\
 &+ \frac{{\left( {y_m^t - y} \right)}}{{\lambda {{\left( {d_{tm}} \right)}^3}}}\left[ {{v_x}\left( {x_m^t - x} \right) + {v_y}\left( {y_m^t - y} \right)} \right]\notag\\
  &+ \frac{{\left( {y_n^r - y} \right)}}{{\lambda {{\left( {d_{rn}} \right)}^3}}}\left[ {{v_x}\left( {x_n^r - x} \right) + {v_y}\left( {y_n^r - y} \right)} \right],
\end{align}
\begin{align}\label{}
{\beta _{nm}} = \frac{{\partial {f_{nm}}}}{{\partial {v_x}}} =
\frac{{x_m^t - x}}{{\lambda d_{tm}}} + \frac{{x_n^r - x}}{{\lambda
d_{rn}}},
\end{align}
\begin{equation}\label{}
{\kappa_{nm}} = \frac{{\partial {f_{nm}}}}{{\partial {v_y}}} =
\frac{{y_m^t - y}}{{\lambda d_{tm}}} + \frac{{y_n^r - y}}{{\lambda
d_{rn}}},
\end{equation}
\begin{equation}\label{}
{\upsilon _{tm}} = \frac{{\partial {d_{tm}}}}{{\partial x}} =
\frac{{x - x_m^t}}{{{d_{tm}}}},
\end{equation}
\begin{equation}\label{}
{l_{tm}} = \frac{{\partial {d_{tm}}}}{{\partial y}} = \frac{{y -
y_m^t}}{{{d_{tm}}}},
\end{equation}
\begin{equation}\label{}
{\eta _{rn}} = \frac{{\partial {d_{rn}}}}{{\partial x}} = \frac{{x -
x_n^r}}{{{d_{rn}}}},
\end{equation}
and
\begin{equation}\label{}
{\psi _{rn}} = \frac{{\partial {d_{rn}}}}{{\partial y}} = \frac{{y -
y_n^r}}{{{d_{rn}}}}.
\end{equation}
 Note that ${a_{nm}}$, ${b_{nm}}$, ${e_{nm}}$,
${g_{nm}}$, ${\beta _{nm}}$, ${\kappa_{nm}}$, $\upsilon _{tm}$,
$l_{tm}$, $\eta _{rn}$ and $\psi _{rn}$ are determined by the target
position and velocity, as well as the position of the receivers and
transmitters.

{\subsection{Calculation of ${ {\bm{J}(\bm{\vartheta}|\bm{\alpha})}
}$}}

According to the likelihood function in (19), the $ij$th element of
the FIM for the parameter vector $\bm{\vartheta}$ is given by
\cite{1993}
\begin{align}
  {\left[ {\bm{J}(\bm{\vartheta}|\bm{\alpha})} \right]_{ij}} =
  \text{Tr}\left( {{{\bm{C}}^{ - 1}}\frac{{\partial {\bm{C}}}}{{\partial {{\vartheta} _i}}}{{\bm{C}}^{ - 1}}\frac{{\partial {\bm{C}}}}{{\partial {{\vartheta} _j}}}} \right)\label{J}.
\end{align}
Using the following identities, \cite{2012}
\begin{align}
{\text{Tr}}\left( {{\bm{ABXY}}} \right) = {\left(
{{\text{vec}}{{\left( {\bm{Y}}^\dag \right)}}} \right)^\dag}\left(
{{{\bm{X}}^\dag} \otimes {\bm{A}}} \right){\text{vec}}\left(
{\bm{B}} \right)
\end{align}
and
\begin{align}
{\text{Tr}}\left( {{\bm{AB}}} \right) = {\text{Tr}}\left(
{{\bm{BA}}} \right),
\end{align}
we can rewrite (\ref{J}) as
\begin{align}
  {\left[ {{\bm{J}}\left( \bm \vartheta |\bm{\alpha} \right)} \right]_{ij}} =&Tr\left( {\frac{{\partial {\bm{C}}}}{{\partial {\vartheta _i}}}{{\bm{C}}^{ - 1}}\frac{{\partial {\bm{C}}}}{{\partial {\vartheta _j}}}{{\bm{C}}^{ - 1}}} \right) \notag\hfill \\
   =& {\left( {\frac{{\partial {{\bm{C}}_{vec}}}}{{\partial {\vartheta _i}}}} \right)^H}\left( {{{\bm{C}}^{ - \dag}} \otimes {{\bm{C}}^{ - 1}}} \right)\left( {\frac{{\partial {{\bm{C}}_{vec}}}}{{\partial {\vartheta _j}}}} \right)\label{Jvarij},
\end{align}
where ${{\bm{C}}_{vec}}= {\text{vec}}\left( {\bm{C}}
\right)\label{Cvec}$. {Calculation of the derivatives and further
simplification of {(\ref{Jvarij})} are provided in Appendix A. Then
we can get the final equation
\begin{align}
{\bm{J}}\left( {{\bm{\theta }}|{\bm{\alpha }}} \right) =
\sum\limits_{p = 1}^N {\sum\limits_{q = 1}^M {\sum\limits_{n = 1}^N
{\sum\limits_{m = 1}^M {\left[ \begin{array}{l}
{A_{11}}{\kern 1pt} {\kern 1pt} {\kern 1pt} {\kern 1pt} {\kern 1pt} {\kern 1pt} {\kern 1pt} {\kern 1pt} {\kern 1pt} {\kern 1pt} {\kern 1pt} {\kern 1pt} {\kern 1pt} {\kern 1pt} {A_{12}}{\kern 1pt} {\kern 1pt} {\kern 1pt} {\kern 1pt} {\kern 1pt} {\kern 1pt} {\kern 1pt} {\kern 1pt} {\kern 1pt} {\kern 1pt} {\kern 1pt} {\kern 1pt} {\kern 1pt} {\kern 1pt} {\kern 1pt} {\kern 1pt} {A_{13}}{\kern 1pt} {\kern 1pt} {\kern 1pt} {\kern 1pt} {\kern 1pt} {\kern 1pt} {\kern 1pt} {\kern 1pt} {\kern 1pt} {\kern 1pt} {\kern 1pt} {\kern 1pt} {\kern 1pt} {A_{14}}\\
{A_{21}}{\kern 1pt} {\kern 1pt} {\kern 1pt} {\kern 1pt} {\kern 1pt} {\kern 1pt} {\kern 1pt} {\kern 1pt} {\kern 1pt} {\kern 1pt} {\kern 1pt} {\kern 1pt} {\kern 1pt} {\kern 1pt} {A_{22}}{\kern 1pt} {\kern 1pt} {\kern 1pt} {\kern 1pt} {\kern 1pt} {\kern 1pt} {\kern 1pt} {\kern 1pt} {\kern 1pt} {\kern 1pt} {\kern 1pt} {\kern 1pt} {\kern 1pt} {\kern 1pt} {\kern 1pt} {A_{23}}{\kern 1pt} {\kern 1pt} {\kern 1pt} {\kern 1pt} {\kern 1pt} {\kern 1pt} {\kern 1pt} {\kern 1pt} {\kern 1pt} {\kern 1pt} {\kern 1pt} {\kern 1pt} {\kern 1pt} A{{\kern 1pt} _{24}}\\
{A_{31}}{\kern 1pt} {\kern 1pt} {\kern 1pt} {\kern 1pt} {\kern 1pt} {\kern 1pt} {\kern 1pt} {\kern 1pt} {\kern 1pt} {\kern 1pt} {\kern 1pt} {\kern 1pt} {\kern 1pt} {\kern 1pt} {A_{32{\kern 1pt} {\kern 1pt} {\kern 1pt} {\kern 1pt} {\kern 1pt} {\kern 1pt} {\kern 1pt} }}{\kern 1pt} {\kern 1pt} {\kern 1pt} {\kern 1pt} {\kern 1pt} {\kern 1pt} {\kern 1pt} {\kern 1pt} {A_{33{\kern 1pt} }}{\kern 1pt} {\kern 1pt} {\kern 1pt} {\kern 1pt} {\kern 1pt} {\kern 1pt} {\kern 1pt} {\kern 1pt} {\kern 1pt} {\kern 1pt} {\kern 1pt} {\kern 1pt} {A_{34}}\\
{A_{41}}{\kern 1pt} {\kern 1pt} {\kern 1pt} {\kern 1pt} {\kern 1pt}
{\kern 1pt} {\kern 1pt} {\kern 1pt} {\kern 1pt} {\kern 1pt} {\kern
1pt} {\kern 1pt} {\kern 1pt} {\kern 1pt} {A_{42}}{\kern 1pt} {\kern
1pt} {\kern 1pt} {\kern 1pt} {\kern 1pt} {\kern 1pt} {\kern 1pt}
{\kern 1pt} {\kern 1pt} {\kern 1pt} {\kern 1pt} {\kern 1pt} {\kern
1pt} {\kern 1pt} {\kern 1pt} {A_{43}}{\kern 1pt} {\kern 1pt} {\kern
1pt} {\kern 1pt} {\kern 1pt} {\kern 1pt} {\kern 1pt} {\kern 1pt}
{\kern 1pt} {\kern 1pt} {\kern 1pt} {\kern 1pt} {\kern 1pt} {A_{44}}
\end{array} \right]} } } },\label{jielun1}
\end{align}
where
\begin{align}
\begin{array}{l}
{{\bm{A}}_{11}} = {a_{pq}}({a_{nm}}{({{\bm{J}}_{{\bm{\tau \tau }}}})_{c,d}} + {e_{nm}}{({{\bm{J}}_{{\bm{f\tau }}}})_{c,d}} + \\
{\kern 1pt} {\kern 1pt} {\kern 1pt} {\kern 1pt} {\kern 1pt} {\kern 1pt} {\kern 1pt} {\kern 1pt} {\kern 1pt} {\kern 1pt} {\kern 1pt} {\kern 1pt} {\kern 1pt} {\kern 1pt} {\kern 1pt} {\kern 1pt} {\kern 1pt} {\kern 1pt} {\kern 1pt} {\kern 1pt} {\kern 1pt} {\kern 1pt} {\kern 1pt} {\kern 1pt} {\kern 1pt} {\kern 1pt} {\kern 1pt} {\kern 1pt} {\upsilon _{tm}}{({{\bm{J}}_{{{\bm{d}}_{\bm{t}}}{\bm{\tau }}}})_{m,d}}/N + {\eta _{rn}}{({{\bm{J}}_{{{\bm{d}}_{\bm{r}}}{\bm{\tau }}}})_{n,d}}/M)\\
{\kern 1pt} {\kern 1pt} {\kern 1pt} {\kern 1pt} {\kern 1pt} {\kern 1pt} {\kern 1pt} {\kern 1pt} {\kern 1pt} {\kern 1pt} {\kern 1pt} {\kern 1pt} {\kern 1pt} {\kern 1pt} {\kern 1pt} {\kern 1pt} {\kern 1pt} {\kern 1pt} {\kern 1pt} {\kern 1pt} {\kern 1pt} {\kern 1pt} {\kern 1pt} {\kern 1pt} {\kern 1pt} {\kern 1pt} {\kern 1pt}  + {e_{pq}}({a_{nm}}{({{\bm{J}}_{{\bm{\tau f}}}})_{c,d}} + {e_{nm}}{({{\bm{J}}_{{\bm{ff}}}})_{c,d}} + \\
{\kern 1pt} {\kern 1pt} {\kern 1pt} {\kern 1pt} {\kern 1pt} {\kern 1pt} {\kern 1pt} {\kern 1pt} {\kern 1pt} {\kern 1pt} {\kern 1pt} {\kern 1pt} {\kern 1pt} {\kern 1pt} {\kern 1pt} {\kern 1pt} {\kern 1pt} {\kern 1pt} {\kern 1pt} {\kern 1pt} {\kern 1pt} {\kern 1pt} {\kern 1pt} {\kern 1pt} {\kern 1pt} {\kern 1pt} {\kern 1pt} {\upsilon _{tm}}{({{\bm{J}}_{{{\bm{d}}_{\bm{t}}}{\bm{f}}}})_{m,d}}/N + {\eta _{rn}}{({{\bm{J}}_{{{\bm{d}}_{\bm{r}}}{\bm{f}}}})_{n,d}}/M)\\
{\kern 1pt} {\kern 1pt} {\kern 1pt} {\kern 1pt} {\kern 1pt} {\kern 1pt} {\kern 1pt} {\kern 1pt} {\kern 1pt} {\kern 1pt} {\kern 1pt} {\kern 1pt} {\kern 1pt} {\kern 1pt} {\kern 1pt} {\kern 1pt} {\kern 1pt} {\kern 1pt} {\kern 1pt} {\kern 1pt} {\kern 1pt} {\kern 1pt} {\kern 1pt} {\kern 1pt} {\kern 1pt}  + {\upsilon _{tq}}/N({a_{nm}}{({{\bm{J}}_{{\bm{\tau }}{{\bm{d}}_{\bm{t}}}}})_{c,q}} + {e_{nm}}{({{\bm{J}}_{{\bm{f}}{{\bm{d}}_{\bm{t}}}}})_{c,q}} + \\
{\kern 1pt} {\kern 1pt} {\kern 1pt} {\kern 1pt} {\kern 1pt} {\kern 1pt} {\kern 1pt} {\kern 1pt} {\kern 1pt} {\kern 1pt} {\kern 1pt} {\kern 1pt} {\kern 1pt} {\kern 1pt} {\kern 1pt} {\kern 1pt} {\kern 1pt} {\kern 1pt} {\kern 1pt} {\kern 1pt} {\kern 1pt} {\kern 1pt} {\kern 1pt} {\kern 1pt} {\kern 1pt} {\upsilon _{tm}}{({{\bm{J}}_{{{\bm{d}}_{\bm{t}}}{{\bm{d}}_{\bm{t}}}}})_{m,q}}/N + {\eta _{rn}}{({{\bm{J}}_{{{\bm{d}}_{\bm{r}}}{{\bm{d}}_{\bm{t}}}}})_{n,q}}/M)\\
{\kern 1pt} {\kern 1pt} {\kern 1pt} {\kern 1pt} {\kern 1pt} {\kern 1pt} {\kern 1pt} {\kern 1pt} {\kern 1pt} {\kern 1pt} {\kern 1pt} {\kern 1pt} {\kern 1pt} {\kern 1pt} {\kern 1pt} {\kern 1pt} {\kern 1pt} {\kern 1pt} {\kern 1pt} {\kern 1pt} {\kern 1pt} {\kern 1pt} {\kern 1pt} {\kern 1pt}  + {\eta _{rp}}/M({a_{nm}}{({{\bm{J}}_{{\bm{\tau }}{{\bm{d}}_{\bm{r}}}}})_{c,p}} + {e_{nm}}{({{\bm{J}}_{{\bm{f}}{{\bm{d}}_{\bm{r}}}}})_{c,p}} + \\
{\kern 1pt} {\kern 1pt} {\kern 1pt} {\kern 1pt} {\kern 1pt} {\kern
1pt} {\kern 1pt} {\kern 1pt} {\kern 1pt} {\kern 1pt} {\kern 1pt}
{\kern 1pt} {\kern 1pt} {\kern 1pt} {\kern 1pt} {\kern 1pt} {\kern
1pt} {\kern 1pt} {\kern 1pt} {\kern 1pt} {\kern 1pt} {\kern 1pt}
{\kern 1pt} {\kern 1pt} {\upsilon
_{tm}}{({{\bm{J}}_{{{\bm{d}}_{\bm{t}}}{{\bm{d}}_{\bm{r}}}}})_{m,p}}/N
+ {\eta
_{rn}}{({{\bm{J}}_{{{\bm{d}}_{\bm{r}}}{{\bm{d}}_{\bm{r}}}}})_{n,p}}/M),
\end{array}
\end{align}
\begin{align}
\begin{array}{l}
{{\bm{A}}_{12}} = {{\bm{A}}_{21}}={b_{pq}}({a_{nm}}{({{\bm{J}}_{{\bm{\tau \tau }}}})_{c,d}} + {e_{nm}}{({{\bm{J}}_{{\bm{f\tau }}}})_{c,d}} + \\
{\kern 1pt} {\kern 1pt} {\kern 1pt} {\kern 1pt} {\kern 1pt} {\kern 1pt} {\kern 1pt} {\kern 1pt} {\kern 1pt} {\kern 1pt} {\kern 1pt} {\kern 1pt} {\kern 1pt} {\kern 1pt} {\kern 1pt} {\kern 1pt} {\kern 1pt} {\kern 1pt} {\kern 1pt} {\kern 1pt} {\kern 1pt} {\kern 1pt} {\kern 1pt} {\kern 1pt} {\kern 1pt} {\kern 1pt} {\kern 1pt} {\kern 1pt} {\upsilon _{tm}}{({{\bm{J}}_{{{\bm{d}}_{\bm{t}}}{\bm{\tau }}}})_{m,d}}/N + {\eta _{rn}}{({{\bm{J}}_{{{\bm{d}}_{\bm{r}}}{\bm{\tau }}}})_{n,d}}/M)\\
{\kern 1pt} {\kern 1pt} {\kern 1pt} {\kern 1pt} {\kern 1pt} {\kern 1pt} {\kern 1pt} {\kern 1pt} {\kern 1pt} {\kern 1pt} {\kern 1pt} {\kern 1pt} {\kern 1pt} {\kern 1pt} {\kern 1pt} {\kern 1pt} {\kern 1pt} {\kern 1pt} {\kern 1pt} {\kern 1pt} {\kern 1pt} {\kern 1pt} {\kern 1pt} {\kern 1pt} {\kern 1pt} {\kern 1pt} {\kern 1pt} {\kern 1pt}  + {g_{pq}}({a_{nm}}{({{\bm{J}}_{{\bm{\tau f}}}})_{c,d}} + {e_{nm}}{({{\bm{J}}_{{\bm{ff}}}})_{c,d}} + \\
{\kern 1pt} {\kern 1pt} {\kern 1pt} {\kern 1pt} {\kern 1pt} {\kern 1pt} {\kern 1pt} {\kern 1pt} {\kern 1pt} {\kern 1pt} {\kern 1pt} {\kern 1pt} {\kern 1pt} {\kern 1pt} {\kern 1pt} {\kern 1pt} {\kern 1pt} {\kern 1pt} {\kern 1pt} {\kern 1pt} {\kern 1pt} {\kern 1pt} {\kern 1pt} {\kern 1pt} {\kern 1pt} {\kern 1pt} {\kern 1pt} {\kern 1pt} {\upsilon _{tm}}{({{\bm{J}}_{{{\bm{d}}_{\bm{t}}}{\bm{f}}}})_{m,d}}/N + {\eta _{rn}}{({{\bm{J}}_{{{\bm{d}}_{\bm{r}}}{\bm{f}}}})_{n,d}}/M)\\
{\kern 1pt} {\kern 1pt} {\kern 1pt} {\kern 1pt} {\kern 1pt} {\kern 1pt} {\kern 1pt} {\kern 1pt} {\kern 1pt} {\kern 1pt} {\kern 1pt} {\kern 1pt} {\kern 1pt} {\kern 1pt} {\kern 1pt} {\kern 1pt} {\kern 1pt} {\kern 1pt} {\kern 1pt} {\kern 1pt} {\kern 1pt} {\kern 1pt} {\kern 1pt} {\kern 1pt} {\kern 1pt} {\kern 1pt} {\kern 1pt} {\kern 1pt} {\kern 1pt}  + {l_{tq}}/N({a_{nm}}{({{\bm{J}}_{{\bm{\tau }}{{\bm{d}}_{\bm{t}}}}})_{c,q}} + {e_{nm}}{({{\bm{J}}_{{\bm{f}}{{\bm{d}}_{\bm{t}}}}})_{c,q}} + \\
{\kern 1pt} {\kern 1pt} {\kern 1pt} {\kern 1pt} {\kern 1pt} {\kern 1pt} {\kern 1pt} {\kern 1pt} {\kern 1pt} {\kern 1pt} {\kern 1pt} {\kern 1pt} {\kern 1pt} {\kern 1pt} {\kern 1pt} {\kern 1pt} {\kern 1pt} {\kern 1pt} {\kern 1pt} {\kern 1pt} {\kern 1pt} {\kern 1pt} {\kern 1pt} {\kern 1pt} {\kern 1pt} {\kern 1pt} {\kern 1pt} {\kern 1pt} {\kern 1pt} {\upsilon _{tm}}{({{\bm{J}}_{{{\bm{d}}_{\bm{t}}}{{\bm{d}}_{\bm{t}}}}})_{m,q}}/N + {\eta _{rn}}{({{\bm{J}}_{{{\bm{d}}_{\bm{r}}}{{\bm{d}}_{\bm{t}}}}})_{n,q}}/M)\\
{\kern 1pt} {\kern 1pt} {\kern 1pt} {\kern 1pt} {\kern 1pt} {\kern 1pt} {\kern 1pt} {\kern 1pt} {\kern 1pt} {\kern 1pt} {\kern 1pt} {\kern 1pt} {\kern 1pt} {\kern 1pt} {\kern 1pt} {\kern 1pt} {\kern 1pt} {\kern 1pt} {\kern 1pt} {\kern 1pt} {\kern 1pt} {\kern 1pt} {\kern 1pt} {\kern 1pt} {\kern 1pt} {\kern 1pt} {\kern 1pt} {\kern 1pt}  + {\psi _{rp}}/M({a_{nm}}{({{\bm{J}}_{{\bm{\tau }}{{\bm{d}}_{\bm{r}}}}})_{c,p}} + {e_{nm}}{({{\bm{J}}_{{\bm{f}}{{\bm{d}}_{\bm{r}}}}})_{c,p}} + \\
{\kern 1pt} {\kern 1pt} {\kern 1pt} {\kern 1pt} {\kern 1pt} {\kern
1pt} {\kern 1pt} {\kern 1pt} {\kern 1pt} {\kern 1pt} {\kern 1pt}
{\kern 1pt} {\kern 1pt} {\kern 1pt} {\kern 1pt} {\kern 1pt} {\kern
1pt} {\kern 1pt} {\kern 1pt} {\kern 1pt} {\kern 1pt} {\kern 1pt}
{\kern 1pt} {\kern 1pt} {\kern 1pt} {\kern 1pt} {\kern 1pt} {\kern
1pt} {\kern 1pt} {\upsilon
_{tm}}{({{\bm{J}}_{{{\bm{d}}_{\bm{t}}}{{\bm{d}}_{\bm{r}}}}})_{m,p}}/N
+ {\eta
_{rn}}{({{\bm{J}}_{{{\bm{d}}_{\bm{r}}}{{\bm{d}}_{\bm{r}}}}})_{n,p}}/M),
\end{array}
\end{align}
\begin{align}
\begin{array}{l}
{A_{13}} = {A_{31}} = {\beta _{pq}}[{a_{nm}}{({{\bm{J}}_{{\bm{\tau f}}}})_{c,d}} + {e_{nm}}{({{\bm{J}}_{{\bm{ff}}}})_{c,d}} + \\
{\kern 1pt} {\kern 1pt} {\kern 1pt} {\kern 1pt} {\kern 1pt} {\kern
1pt} {\kern 1pt} {\kern 1pt} {\kern 1pt} {\kern 1pt} {\kern 1pt}
{\kern 1pt} {\kern 1pt} {\kern 1pt} {\kern 1pt} {\kern 1pt} {\kern
1pt} {\kern 1pt} {\kern 1pt} {\kern 1pt} {\kern 1pt} {\kern 1pt}
{\kern 1pt} {\kern 1pt} {\kern 1pt} {\kern 1pt} {\kern 1pt} {\kern
1pt} {\kern 1pt} {\kern 1pt} {\kern 1pt} {\kern 1pt} {\kern 1pt}
{\kern 1pt} {\kern 1pt} {\kern 1pt} {\kern 1pt} {\kern 1pt} {\kern
1pt} {\kern 1pt} {\kern 1pt} {\kern 1pt} {\kern 1pt} {\kern 1pt}
{\kern 1pt} {\kern 1pt} {\kern 1pt} {\kern 1pt} {\kern 1pt} {\kern
1pt} {\kern 1pt} {\kern 1pt} {\kern 1pt} {\kern 1pt} {\kern 1pt}
{\upsilon _{tm}}{({{\bm{J}}_{{{\bm{d}}_{\bm{t}}}{\bm{f}}}})_{m,d}}/N
+ {\eta _{rn}}{({{\bm{J}}_{{{\bm{d}}_{\bm{r}}}{\bm{f}}}})_{n,d}}/M],
\end{array}
\end{align}
\begin{align}
\begin{array}{l}
{A_{14}} = {A_{41}} = {k_{pq}}[{a_{nm}}{({{\bm{J}}_{{\bm{\tau f}}}})_{c,d}} + {e_{nm}}{({{\bm{J}}_{{\bm{ff}}}})_{c,d}} + \\
{\kern 1pt} {\kern 1pt} {\kern 1pt} {\kern 1pt} {\kern 1pt} {\kern
1pt} {\kern 1pt} {\kern 1pt} {\kern 1pt} {\kern 1pt} {\kern 1pt}
{\kern 1pt} {\kern 1pt} {\kern 1pt} {\kern 1pt} {\kern 1pt} {\kern
1pt} {\kern 1pt} {\kern 1pt} {\kern 1pt} {\kern 1pt} {\kern 1pt}
{\kern 1pt} {\kern 1pt} {\kern 1pt} {\kern 1pt} {\kern 1pt} {\kern
1pt} {\kern 1pt} {\kern 1pt} {\kern 1pt} {\kern 1pt} {\kern 1pt}
{\kern 1pt} {\kern 1pt} {\kern 1pt} {\kern 1pt} {\kern 1pt} {\kern
1pt} {\kern 1pt} {\kern 1pt} {\kern 1pt} {\kern 1pt} {\kern 1pt}
{\kern 1pt} {\kern 1pt} {\kern 1pt} {\kern 1pt} {\kern 1pt} {\kern
1pt} {\kern 1pt} {\kern 1pt} {\kern 1pt} {\kern 1pt} {\kern 1pt}
{\upsilon _{tm}}{({{\bm{J}}_{{{\bm{d}}_{\bm{t}}}{\bm{f}}}})_{m,d}}/N
+ {\eta _{rn}}{({{\bm{J}}_{{{\bm{d}}_{\bm{r}}}{\bm{f}}}})_{n,d}}/M],
\end{array}
\end{align}
\begin{align}
\begin{array}{l}
{{A}_{22}} = {b_{pq}}({b_{nm}}{({{\bm{J}}_{{\bm{\tau \tau }}}})_{c,d}} + {g_{nm}}{({{\bm{J}}_{{\bm{f\tau }}}})_{c,d}} + \\
{\kern 1pt} {\kern 1pt} {\kern 1pt} {\kern 1pt} {\kern 1pt} {\kern 1pt} {\kern 1pt} {\kern 1pt} {\kern 1pt} {\kern 1pt} {\kern 1pt} {\kern 1pt} {\kern 1pt} {\kern 1pt} {\kern 1pt} {\kern 1pt} {\kern 1pt} {\kern 1pt} {\kern 1pt} {\kern 1pt} {\kern 1pt} {\kern 1pt} {\kern 1pt} {\kern 1pt} {\kern 1pt} {\kern 1pt} {\kern 1pt} {\kern 1pt} {l_{tm}}{({{\bm{J}}_{{{\bm{d}}_{\bm{t}}}{\bm{\tau }}}})_{m,d}}/N + {\psi _{rn}}{({{\bm{J}}_{{{\bm{d}}_{\bm{r}}}{\bm{\tau }}}})_{n,d}}/M)\\
{\kern 1pt} {\kern 1pt} {\kern 1pt} {\kern 1pt} {\kern 1pt} {\kern 1pt} {\kern 1pt} {\kern 1pt} {\kern 1pt} {\kern 1pt} {\kern 1pt} {\kern 1pt} {\kern 1pt} {\kern 1pt} {\kern 1pt} {\kern 1pt} {\kern 1pt} {\kern 1pt} {\kern 1pt} {\kern 1pt} {\kern 1pt} {\kern 1pt} {\kern 1pt} {\kern 1pt} {\kern 1pt} {\kern 1pt} {\kern 1pt} {\kern 1pt}  + {g_{pq}}({b_{nm}}{({{\bm{J}}_{{\bm{\tau f}}}})_{c,d}} + {g_{nm}}{({{\bm{J}}_{{\bm{ff}}}})_{c,d}} + \\
{\kern 1pt} {\kern 1pt} {\kern 1pt} {\kern 1pt} {\kern 1pt} {\kern 1pt} {\kern 1pt} {\kern 1pt} {\kern 1pt} {\kern 1pt} {\kern 1pt} {\kern 1pt} {\kern 1pt} {\kern 1pt} {\kern 1pt} {\kern 1pt} {\kern 1pt} {\kern 1pt} {\kern 1pt} {\kern 1pt} {\kern 1pt} {\kern 1pt} {\kern 1pt} {\kern 1pt} {\kern 1pt} {\kern 1pt} {\kern 1pt} {l_{tm}}{({{\bm{J}}_{{{\bm{d}}_{\bm{t}}}{\bm{f}}}})_{m,d}}/N + {\psi _{rn}}{({{\bm{J}}_{{{\bm{d}}_{\bm{r}}}{\bm{f}}}})_{n,d}}/M)\\
{\kern 1pt} {\kern 1pt} {\kern 1pt} {\kern 1pt} {\kern 1pt} {\kern 1pt} {\kern 1pt} {\kern 1pt} {\kern 1pt} {\kern 1pt} {\kern 1pt} {\kern 1pt} {\kern 1pt} {\kern 1pt} {\kern 1pt} {\kern 1pt} {\kern 1pt} {\kern 1pt} {\kern 1pt} {\kern 1pt} {\kern 1pt} {\kern 1pt} {\kern 1pt} {\kern 1pt} {\kern 1pt} {\kern 1pt}  + {l_{tq}}/N({b_{nm}}{({{\bm{J}}_{{\bm{\tau }}{{\bm{d}}_{\bm{t}}}}})_{c,q}} + {g_{nm}}{({{\bm{J}}_{{\bm{f}}{{\bm{d}}_{\bm{t}}}}})_{c,q}} + \\
{\kern 1pt} {\kern 1pt} {\kern 1pt} {\kern 1pt} {\kern 1pt} {\kern 1pt} {\kern 1pt} {\kern 1pt} {\kern 1pt} {\kern 1pt} {\kern 1pt} {\kern 1pt} {\kern 1pt} {\kern 1pt} {\kern 1pt} {\kern 1pt} {\kern 1pt} {\kern 1pt} {\kern 1pt} {\kern 1pt} {\kern 1pt} {\kern 1pt} {\kern 1pt} {\kern 1pt} {\kern 1pt} {\kern 1pt} {l_{tm}}{({{\bm{J}}_{{{\bm{d}}_{\bm{t}}}{{\bm{d}}_{\bm{t}}}}})_{m,q}}/N + {\psi _{rn}}{({{\bm{J}}_{{{\bm{d}}_{\bm{r}}}{{\bm{d}}_{\bm{t}}}}})_{n,q}}/M)\\
{\kern 1pt} {\kern 1pt} {\kern 1pt} {\kern 1pt} {\kern 1pt} {\kern 1pt} {\kern 1pt} {\kern 1pt} {\kern 1pt} {\kern 1pt} {\kern 1pt} {\kern 1pt} {\kern 1pt} {\kern 1pt} {\kern 1pt} {\kern 1pt} {\kern 1pt} {\kern 1pt} {\kern 1pt} {\kern 1pt} {\kern 1pt} {\kern 1pt} {\kern 1pt} {\kern 1pt} {\kern 1pt}  + {\psi _{rp}}/M({b_{nm}}{({{\bm{J}}_{{\bm{\tau }}{{\bm{d}}_{\bm{r}}}}})_{c,p}} + {g_{nm}}{({{\bm{J}}_{{\bm{f}}{{\bm{d}}_{\bm{r}}}}})_{c,p}} + \\
{\kern 1pt} {\kern 1pt} {\kern 1pt} {\kern 1pt} {\kern 1pt} {\kern
1pt} {\kern 1pt} {\kern 1pt} {\kern 1pt} {\kern 1pt} {\kern 1pt}
{\kern 1pt} {\kern 1pt} {\kern 1pt} {\kern 1pt} {\kern 1pt} {\kern
1pt} {\kern 1pt} {\kern 1pt} {\kern 1pt} {\kern 1pt} {\kern 1pt}
{\kern 1pt} {\kern 1pt} {\kern 1pt}
{l_{tm}}{({{\bm{J}}_{{{\bm{d}}_{\bm{t}}}{{\bm{d}}_{\bm{r}}}}})_{m,p}}/N
+ {\psi
_{rn}}{({{\bm{J}}_{{{\bm{d}}_{\bm{r}}}{{\bm{d}}_{\bm{r}}}}})_{n,p}}/M)
\end{array}
\end{align}
\begin{align}
\begin{array}{l}
{A_{23}} = {A_{32}} = {\beta _{pq}}[{b_{nm}}{({{\bm{J}}_{{\bm{\tau f}}}})_{c,d}} + {g_{nm}}{({{\bm{J}}_{{\bm{ff}}}})_{c,d}} + \\
{\kern 1pt} {\kern 1pt} {\kern 1pt} {\kern 1pt} {\kern 1pt} {\kern
1pt} {\kern 1pt} {\kern 1pt} {\kern 1pt} {\kern 1pt} {\kern 1pt}
{\kern 1pt} {\kern 1pt} {\kern 1pt} {\kern 1pt} {\kern 1pt} {\kern
1pt} {\kern 1pt} {\kern 1pt} {\kern 1pt} {\kern 1pt} {\kern 1pt}
{\kern 1pt} {\kern 1pt} {\kern 1pt} {\kern 1pt} {\kern 1pt} {\kern
1pt} {\kern 1pt} {\kern 1pt} {\kern 1pt} {\kern 1pt} {\kern 1pt}
{\kern 1pt} {\kern 1pt} {\kern 1pt} {\kern 1pt} {\kern 1pt} {\kern
1pt} {\kern 1pt} {\kern 1pt} {\kern 1pt} {\kern 1pt} {\kern 1pt}
{\kern 1pt} {\kern 1pt} {\kern 1pt} {\kern 1pt} {\kern 1pt} {\kern
1pt} {\kern 1pt} {\kern 1pt} {\kern 1pt} {\kern 1pt}
{l_{tm}}{({{\bm{J}}_{{{\bm{d}}_{\bm{t}}}{\bm{f}}}})_{m,d}}/N + {\psi
_{rn}}{({{\bm{J}}_{{{\bm{d}}_{\bm{r}}}{\bm{f}}}})_{n,d}}/M],
\end{array}
\end{align}
\begin{align}
\begin{array}{l}
{A_{24}} = {A_{42}} = {k_{pq}}[{b_{nm}}{({{\bm{J}}_{{\bm{\tau f}}}})_{c,d}} + {g_{nm}}{({{\bm{J}}_{{\bm{ff}}}})_{c,d}} + \\
{\kern 1pt} {\kern 1pt} {\kern 1pt} {\kern 1pt} {\kern 1pt} {\kern
1pt} {\kern 1pt} {\kern 1pt} {\kern 1pt} {\kern 1pt} {\kern 1pt}
{\kern 1pt} {\kern 1pt} {\kern 1pt} {\kern 1pt} {\kern 1pt} {\kern
1pt} {\kern 1pt} {\kern 1pt} {\kern 1pt} {\kern 1pt} {\kern 1pt}
{\kern 1pt} {\kern 1pt} {\kern 1pt} {\kern 1pt} {\kern 1pt} {\kern
1pt} {\kern 1pt} {\kern 1pt} {\kern 1pt} {\kern 1pt} {\kern 1pt}
{\kern 1pt} {\kern 1pt} {\kern 1pt} {\kern 1pt} {\kern 1pt} {\kern
1pt} {\kern 1pt} {\kern 1pt} {\kern 1pt} {\kern 1pt} {\kern 1pt}
{\kern 1pt} {\kern 1pt} {\kern 1pt} {\kern 1pt} {\kern 1pt} {\kern
1pt} {\kern 1pt} {\kern 1pt} {\kern 1pt} {\kern 1pt}
{l_{tm}}{({{\bm{J}}_{{{\bm{d}}_{\bm{t}}}{\bm{f}}}})_{m,d}}/N + {\psi
_{rn}}{({{\bm{J}}_{{{\bm{d}}_{\bm{r}}}{\bm{f}}}})_{n,d}}/M],
\end{array}
\end{align}
\begin{align}
{{{A}}_{33}} ={\beta _{pq}}{\beta
_{nm}}{({{\bm{J}}_{{\bm{ff}}}})_{c,d}},
\end{align}
\begin{align}
{{{A}}_{34}} ={{{A}}_{43}} = {k_{pq}}{\beta
_{nm}}{({{\bm{J}}_{{\bm{ff}}}})_{c,d}},
\end{align}
and
\begin{align}
{{{A}}_{44}} ={\kern 1pt}
{k_{pq}}{k_{nm}}{({{\bm{J}}_{{\bm{ff}}}})_{c,d}},\label{jielun2}
\end{align}
{where $c = M(n - 1) + m$ and $d = M(p - 1) + q$.
${{\bm{J}}_{{\bm{\tau \tau }}}}$, ${{\bm{J}}_{{\bm{\tau f}}}}$,
${{\bm{J}}_{{\bm{f \tau}}}}$, ${{\bm{J}}_{{\bm{\tau
}}{{\bm{d}}_{\bm{t}}}}}$, ${{\bm{J}}_{{{\bm{d}}_{\bm{t}}}{\bm{\tau
}}}}$, ${{\bm{J}}_{{\bm{\tau }}{{\bm{d}}_{\bm{r}}}}}$,
${{\bm{J}}_{{{\bm{d}}_{\bm{r}}}{\bm{\tau }}}}$,
${{\bm{J}}_{{\bm{ff}}}}$,
${{\bm{J}}_{{\bm{f}}{{\bm{d}}_{\bm{t}}}}}$,
${{\bm{J}}_{{{\bm{d}}_{\bm{t}}}{\bm{f}}}}$,
${{\bm{J}}_{{\bm{f}}{{\bm{d}}_{\bm{r}}}}}$,
${{\bm{J}}_{{{\bm{d}}_{\bm{r}}}{\bm{f}}}}$,
${{\bm{J}}_{{{\bm{d}}_{\bm{t}}}{{\bm{d}}_{\bm{t}}}}}$,
${{\bm{J}}_{{{\bm{d}}_{\bm{t}}}{{\bm{d}}_{\bm{r}}}}}$,
${{\bm{J}}_{{{\bm{d}}_{\bm{r}}}{{\bm{d}}_{\bm{t}}}}}$,
${{\bm{J}}_{{{\bm{d}}_{\bm{r}}}{{\bm{d}}_{\bm{r}}}}}$ are defined in
(72).} {It should be noted that, the results obtained here, say
(44)-(54), are a highly nontrival extension of the previous results
in \cite{He:2010}. Unfortunately, they are, as one might expect,
considerably more complicated but they describe the best possible
estimation performance in non-ideal scenarios that are of great
practical interest in the following sense.} Given any unbiased
estimator $ \hat{ {\boldsymbol \theta } } $ of an unknown parameter
$ {\boldsymbol \theta} $ based on an observation vector $ {\bm r} $,
when ${\bm \alpha }$ is assumed known and fixed, we have \cite{1993}
\begin{equation}
\textup{MSE}={\mathbb{E}_{{\bm{r}}|{\boldsymbol{\theta
}},{\boldsymbol{\alpha }}}}\left\{ ( \hat{ {\boldsymbol \theta} } -
{\bm{\theta }} ) ( \hat{ {\boldsymbol \theta} } - {\bm{\theta
}})^\dag
 \right\} \succeq
 \textup{CRB}({{\bm{\theta }}}|{\bm{\alpha }}) =
{{\bm{J}}^{ - 1}}({{\bm{\theta }}}|{\bm{\alpha }}).
\label{crb-bound}
\end{equation}
which is the standard CRB for vector parameters where $ {\bm A }
\succeq {\bm B } $ means $ {\bm A } - {\bm B } $ is positive
semidefinite, and MSE is the mean squared error {matrix} of the
unbiased estimator.

\section{Cramer-Rao Bound For Mismatched Case}
\label{mm}

In order to find an ML estimate or use the CRB result in (55), now
called the generalized CRB (GCRB), we must know the actual values of
the signal matrix $ {\bm S} $ (and thus ${\bm \alpha}$) from (13),
the reflection coefficients covariance matrix ${\bm R }$ described
near (14), and the noise covariance matrix ${\bm Q }$ described near
(15). Here, we assume the estimation algorithm employs incorrect
values for these matrices denoted by ${{\bm{S}}_0}$, ${{\bm{R}}_0}$,
{and} ${{\bm{Q}}_0}$ respectively. The incorrect values
${{\bm{S}}_0}$, ${{\bm{R}}_0}$, {and} ${{\bm{Q}}_0}$ might be
obtained from some inaccurate estimation. Given the estimation
algorithm uses these incorrect values ${{\bm{S}}_0}$,
${{\bm{R}}_0}$, {and} ${{\bm{Q}}_0}$, we find a lower bound on the
estimation performance using some recently published work
\cite{mismatch}. In the case described, the assumed likelihood
function is
\begin{align}
{p_0}({\bm{r}}|{\bm{\theta }},{\bm{\alpha }}) = \frac{1}{{{\pi
^{KN}}\det {\bm{C_0}}}}\exp ( - {{\bm{r}}^H}{{\bm{C_0}}^{ -
1}}{\bm{r}}),
\end{align}
where ${{\bm{C}}_0} = {{\bm{S}}_0}{{\bm{R}}_0}{{\bm{S}}_0}^H +
{{\bm{Q}}_0}$. To avoid confusion with the GCRB, we {denote} the
actual values of the signal matrix from (13), the reflection
coefficients covariance matrix described near (14), and noise
covariance matrix described near (15) {by} ${{\bm{S}}_1}$,
${{\bm{R}}_1}$, {and} ${{\bm{Q}}_1}$.

Thus, the actual likelihood function is
\begin{align}
{p_1}({\bm{r}}|{\bm{\theta }},{\bm{\alpha }}) = \frac{1}{{{\pi
^{KN}}\det {{\bm{C}}_1}}}\exp ( - {{\bm{r}}^H}{{\bm{C}}_1}^{ -
1}{\bm{r}})
\end{align}
where ${{\bm{C}}_1} = {{\bm{S}}_1}{{\bm{R}}_1}{{\bm{S}}_1}^H +
{{\bm{Q}}_1}$. According to \cite{mismatch}, we know that
\begin{equation}
\textup{MSE}_{mis}\succeq
  \textup{CRB}_{mis}({{\bm{\theta }}}|{\bm{\alpha }}) =
{{\bm{J_{mis}}}^{ - 1}}({{\bm{\theta }}}|{\bm{\alpha }}).
\end{equation}
{where $\rm{MSE}_{mis}$, $\rm{CR{B_{mis}}({{\bm{\theta
}}}|{\bm{\alpha }})}$ and ${{\bm{J_{mis}}}}({{\bm{\theta
}}}|{\bm{\alpha }})$ {denote }the MSE, CRB and FIM {matrices} under
mismatched situation, and }
\begin{align}
{{\bm{J}}_{mis}}({\bm{\theta }}|{\bm{\alpha }}) =
&{\mathbb{E}_{{p_1}({\bm{r}}|{\bm{\theta }},{\bm{\alpha }})}}
\biggl{\{} { \biggl{(} \frac{{{p_0}({\bm{r}}|{\bm{\theta }},
{\bm{\alpha }})}}{{{p_1}({\bm{r}}|{\bm{\theta }},{\bm{\alpha }})}} \biggr{)}^2}\notag\\
&\times {\nabla _{\bm{\theta }}}\log {p_0}({\bm{r}}|{\bm{\theta
}},{\bm{\alpha }}){[{\nabla _{\bm{\theta }}}\log
{p_0}({\bm{r}}|{\bm{\theta }},{\bm{\alpha }})]^\dag } \biggr{\}}
\end{align}

Next note that ${\nabla _{\bm{\theta }}}\log p({\bm{r}}|{\bm{\theta
}},{\bm{\alpha }}) = ({\nabla _\theta }{{\bm{\vartheta}} ^\dag }){\nabla
_{\bm{\vartheta}} }\log p({\bm{r}}|{\bm{\vartheta}} ,{\bm{\alpha }})$ and
${p_0}({\bm{r}}|{\bm{\theta }},{\bm{\alpha }}) =
{p_0}({\bm{r}}|{\bm{\vartheta}} ,{\bm{\alpha }})$ so that
\begin{align}
{{\bm{J}}_{mis}}({\bm{\theta }}|{\bm{\alpha }}) = ({\nabla
_{\bm{\theta }}}{{\bm{\vartheta}} ^\dag }){{\bm{J}}_{mis}}({\bm{\vartheta}}
|{\bm{\alpha }}){({\nabla _{\bm{\theta }}}{{\bm{\vartheta}} ^\dag })^\dag
},
\end{align}
where
\begin{align}
{{\bm{J}}_{mis}}({\bm{\vartheta}} |{\bm{\alpha }}) = &{\mathbb{E}_{{p_1}({\bm{r}}|{\bm{\vartheta}} ,{\bm{\alpha }})}} \biggl{\{} { \biggl{(} \frac{{{p_0}({\bm{r}}|{\bm{\vartheta}} ,{\bm{\alpha }})}}{{{p_1}({\bm{r}}|{\bm{\vartheta}} ,{\bm{\alpha }})}} \biggr{)}^2}\notag\\
&\times{\nabla _{\bm{\vartheta}} }\log {p_0}({\bm{r}}|{\bm{\vartheta}} ,{\bm{\alpha
}}){[{\nabla _{\bm{\vartheta}} }\log {p_0}({\bm{r}}|{\bm{\vartheta}} ,{\bm{\alpha
}})]^\dag } \biggr{\}} \label{a}
\end{align}

Calculation of the derivatives and further simplification of
(\ref{a}) is omitted due to similarity to the case without mismatch.

\section{Numerical Examples}\label{sec:Result}

In this section, examples are presented which demonstrate the use of
the GCRB {and the mismatched CRB} presented in the previous section
to bound the performance of distributed radar networks which employ
multiple widely spaced transmitters and receivers to jointly
estimate target position and velocity.  For brevity, we focus on
examples which employ signals that are more applicable for passive
radar. Initially, we describe performance when the transmitted
signals are either known or where the transmitted signals of
opportunity are estimated perfectly from a direct path reception.
Later we consider cases where this is not true.  We also assume that
the positions of the transmitters and receivers are exactly known.
For passive radar cases, these assumptions allow us to describe the
best possible performance that can be obtained under the best
circumstances.

It is easy to employ our bounds for cases where all parameters to be
included in the vector $ {\bm{\alpha }} $ are known and thus the bound in
(\ref{crb-bound}) is applicable.  However, the vector $ {\bm{\alpha }} $
might include a random bit sequence which contains information being transmitted.
In order to avoid presenting a CRB for every possible bit sequence (${\bf \alpha} $),
we quote the expected CRB averaged over all bit sequences (ECRBOB), assuming each
bit sequence to be perfectly estimated. From (\ref{crb-bound}),
\begin{align}
ECRBOB({{\bm{\theta }}}) = {\mathbb{E}}_{\bm{\alpha }}\left\{
{CR{B}(\bm{\theta }\left| \bm{\alpha }\right.) } \right\}
\end{align}
clearly bounds the corresponding covariance matrrix averaged over
all bit sequences. For the best case, when the bit sequence in $ {\bm{\alpha }} $ is
perfectly estimated, the ECRBOB is a good indicator of performance.
For example, it describes how the system parameters, such as the
number of antennas, the geometry, and the waveforms impact
performance, assuming accurate estimation of $\bm{\alpha}$. One can
use the ECRBOB to optimize any parameter of interest.

Consider a target moving with velocity $(50,30)$ m/s is present at
$( 15.15,10.1275)$ km. To define a general test set up that is easy
to describe for general $M$ and $N$, each transmit and receive
(single antenna) station is located 7 km from the reference point
(15,10) km. The $M$ transmit stations are uniformly distributed in
angle over the range $[0,2\pi)$, i.e, the angle of the $m$-th
{transmitter} is $\varphi _m^t = 2\pi (m - 1)/M, m = 1, \cdots ,M$.
The $N$ receive stations are also uniformly distributed in angle
over the range $[0,2\pi)$, i.e. the angle of the $n$-th receive
station is $\varphi_n^r = 2\pi (n - 1)/N, n = 1, \cdots ,N$, where
the angles are measured with respect to the horizontal axis
originated at the reference point as illustrated in Figure 1.
Suppose $E_1=E_2=...=E_M=E$. Fix $SCNR = 10\log ((\sum\nolimits_{n =
1}^N {\sum\nolimits_{m = 1}^M {\sigma
_{nm}^2E{P_0}/d_{tm}^2d_{rn}^2)} /(N\sigma _w^2)} $, called the
signal-to-clutter-plus-noise ratio (SCNR), where ${{\sigma_{nm}^2}}
= \mathbb{E}\{{\zeta_{nm}}{\zeta_{nm}}^H \}$ and ${{\sigma_w^2}} =
\mathbb{E}\{{w_n(1)}{w_n(1)}^H
\}=\cdots=\mathbb{E}\{{w_n(K)}{w_n(K)}^H \}$. We set
$\sigma_{nm}^2=1$ for all $n$ and $m$, and ${P_0} = 1$.

\begin{figure}
\centering
\includegraphics[width=3.0 in]{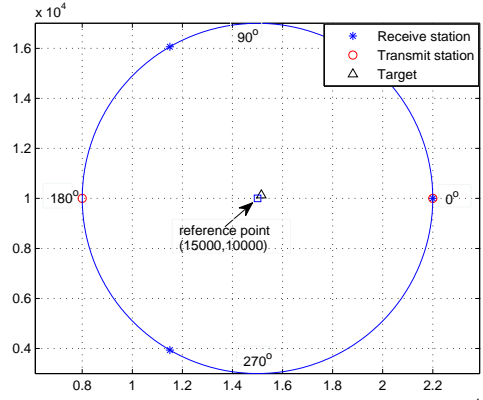}
\vspace{-0.4cm} \caption{Parameter set up for a distributed radar
network with $M=2$ and $N=3$.}\label{MIMO} \vspace{-0.1cm}
\end{figure}

To be relevant to a passive radar system, the signals considered are
those employed by the popular Global System for Mobile (GSM)
Communications system. The baseband transmitted waveforms are
Gaussian minimum shift keying (GMSK) signals [13]
\begin{align}
{s_m}(k,{{\bm{\alpha }}_m}) = {A_m}\exp \left\{ {j\sum\limits_{i =
1}^{{N_c}} {{c_{mi}}} \sum\limits_{j = 1}^k {z(k{T_s} -
i{T_p}){T_s}} } \right\}{e^{j2\pi m \triangle fkT_s}},
\label{GSM}
\end{align}
where ${{\bm{\alpha }}_m} = {\left[ {{c_{m1}}, \ldots
,{c_{m{N_{_c}}}}} \right]^\dag }$,
\begin{align}
z(t) = \frac{\pi }{{2{T_p}}}\left\{ {{{\vartheta}} \left[ {\frac{{2\pi
B}}{{\sqrt {\ln 2} }}(t - {T_P})} \right] - {{\vartheta}} \left[
{\frac{{2\pi B}}{{\sqrt {\ln 2} }}}t \right]} \right\},
\end{align}
$\vartheta [t] = \left( {1/\sqrt {2\pi } } \right)\int_t^\infty
{{e^{ - {\tau ^2}/2}}} d\tau $, 
$T_p$ is the bit duration, ${B}$ denotes the 3 dB bandwidth of the
Gaussian prefilter used in the GMSK modulators, ${c_{mi}} \in \{ -
1,1\}$ is the $i$th $(i = 1, \ldots ,{N_c})$ binary data bit of the
$m$th transmitted waveform,  ${N_c}$
denotes the number of bits contained in the observation interval,
${A_m}$ is the normalization factor, and $\triangle f = f_{k+1}-f_k$
is the frequency offset between different signals of opportunity
with neighboring frequencies.  In the simulations, we generate
${c_{mi}}=-1$ or 1 randomly with the same probability of $0.5$.
To model a GSM system, assume $
T_p=577\mu s, BT_p=0.3, N_c=16 $, the carrier frequency $f_c=900$
MHz and $\triangle f=3$ KHz (orthogonal signals) or $\triangle
f=300$ Hz (nonorthogonal signals). {It should be noticed that the
bandwidth is only 520Hz and $N_c=16$ in the simulation because of
the huge calculation complexity.}

\begin{figure}
\centering
 \includegraphics[width=3.0 in]{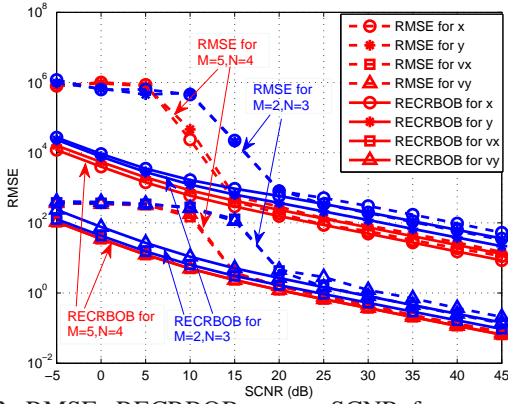}
\vspace{-0.4cm} \caption{{RMSE, RECRBOB versus SCNR for a passive
distributed radar network with $M=2$, $N=3$ and $M=5$, $N=4$,
spatially independent refection coefficients, spatially independent
noise, and nonorthogonal signal}.} \vspace{-0.1cm}\label{5423}
\end{figure}

{Figure \ref{5423} shows the cases with $M=2$, $N=3$ and $M=5$,
$N=4$ for spatially independent refection coefficients, spatially
independent noise, and nonorthogonal signals. The solid and dashed
curves show the root ECRBOB (RECRBOB) and the root mean squared
error (RMSE) of the ML estimation, respectively, in the cases
investigated. It is seen that all curves show that the RMSE
decreases as the signal-to-clutter-plus-noise ratio (SCNR) is
increased.  In support of the correctness of our derived CRBs, all
RMSE curves show the existence of a threshold, above which the RMSE
starts to become close to the RECRBOB in value and slope. We can see
that the threshold in the case with $M=5$, $N=4$ is 15 dB, while the
threshold in the case with $M=2$, $N=3$ is 20 dB. The MSE curves for
$M=5$, $N=4$ have a significantly lower threshold than those for
$M=2$, $N=3$, apparently due to the additional transmit and receive
stations, while the reduction in RECRBOB and the RMSE above
threshold due to employing $M=5$, $N=4$ instead of $M=2$, $N=3$ is
significantly smaller. Also, in the high SCNR region, RMSE is closer
to RECRBOB for the case with $M=5$, $N=4$ than for the case with
$M=2$, $N=3$.}

Increasing the time duration of the signals, $N_c$ in (\ref{GSM})
also provides benefits as one would expect. While Figure \ref{5423}
considers the case of $N_c=16$, Figure \ref{5} shows a comparison
between RECRBOB for $N_c=16$ and  RECRBOB for $N_c=64$ for the same
case of $M=2$, $N=3$, spatially independent refection coefficients,
spatially independent noise, and orthogonal signals.  Here we can
see that if we increase $N_c$, the RECRBOB will be significant reduced.
In the following cases, to reduce complexity, we employ $N_c=16$.

\begin{figure}
\centering
 \includegraphics[width=3.0 in]{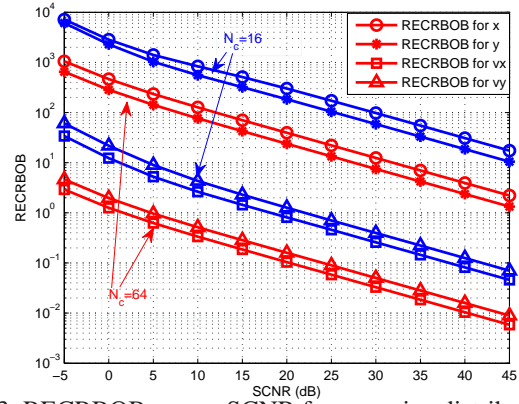}
\vspace{-0.4cm} \caption{{{RECRBOB versus SCNR for a passive
distributed radar network with $M=2$, $N=3$, spatially independent
refection coefficients, spatially independent noise, and orthogonal
signals}}.} \vspace{-0.1cm}\label{5}
\end{figure}

\subsubsection{{Orthogonal Signals and Nonorthogonal Signals}}

In this section, we focus on the effect of the nonorthogonality of the
different transmitted signals. We consider the situation of
spatially independent reflection coefficients, spatially independent
noise, and nonorthogonal signals.  The other factors are the same as in
Figure \ref{5423}. The system considered in Figure
\ref{orthigonality} has $M=2$ transmitters and $N=3$ receivers. The
red and blue curves in this figure correspond to the cases with
orthogonal and nonorthogonal signals, respectively. We see that the
threshold obtained using the orthogonal signals is 15 dB while the
threshold obtained using the nonorthogonal signals is 20 dB. Thus
the threshold for the nonorthogonal signals tested is higher than
the threshold for the orthogonal signals tested. It is also seen
that the RECRBOB of the orthogonal signals tested is smaller than
that for the nonorthogonal signals tested over the whole region of
SCNR shown. So both the RMSE and RECRBOB indicate that the radar can
achieve better performance if the waveforms are closer to being
orthogonal in the case considered.

\begin{figure}
\centering
\includegraphics[width=3.0 in]{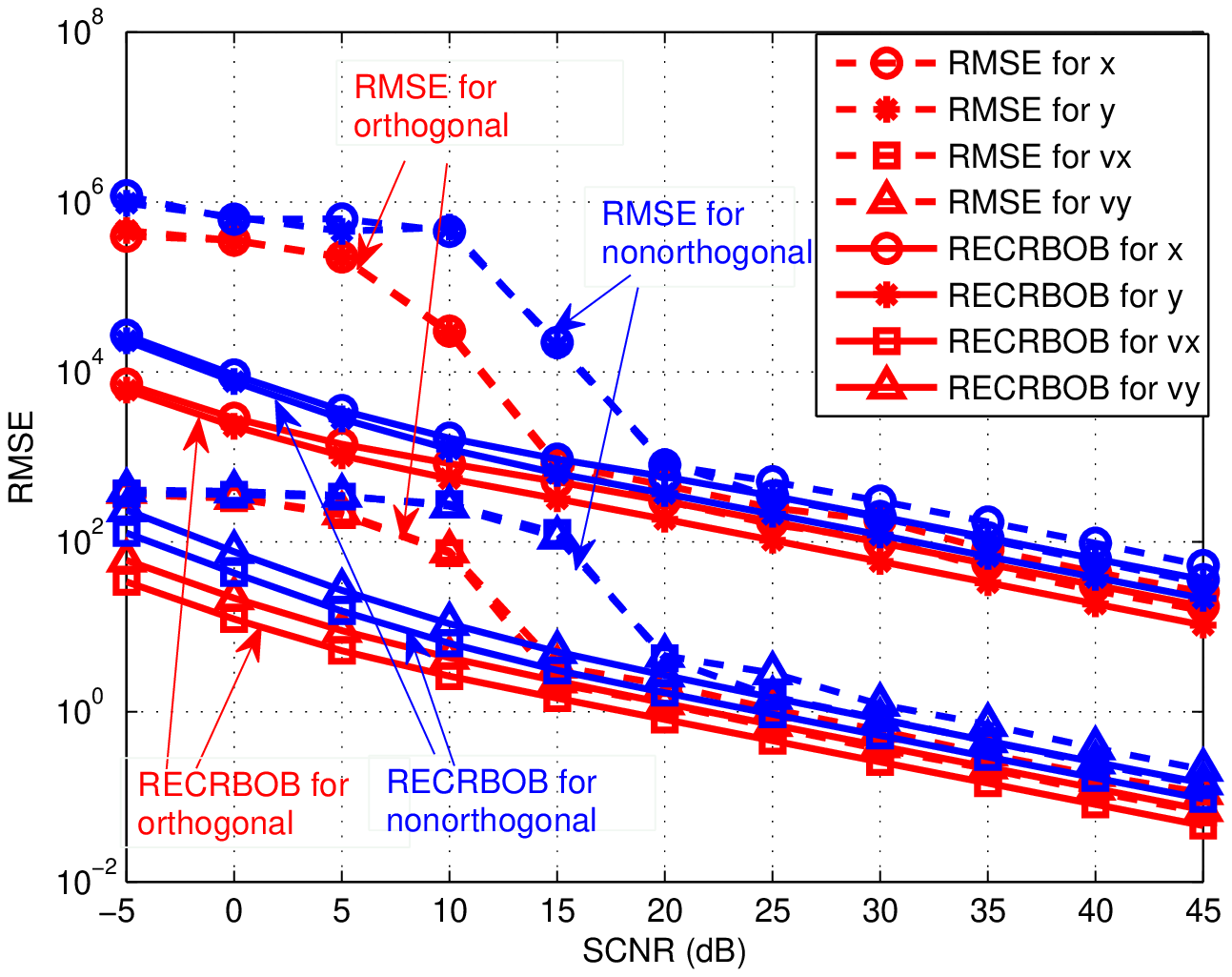}
\vspace{-0.4cm} \caption{RMSE, RECRBOB versus SCNR for a passive
MIMO radar with $M=2$ and $N=3$, spatially independent refection
coefficients, and spatially independent noise.}\label{orthigonality}
\vspace{-0.1cm}
\end{figure}

\subsubsection{Spatially Dependent Reflection Coefficients}

In this section, we consider the situation of spatially dependent
reflection coefficients, spatially independent noise, and orthogonal
signals. The elements of the covariance matrix $\bm{R}$ describing
the correlation between the different reflection coefficients are
generated with {\cite{He:2010}}
\begin{align}
\bm{R}={\bm{R}^r}\otimes{\bm{R}^t}
\end{align}
where
\begin{align}
{{\bm{R}}^r} = \left[ {\begin{array}{*{20}{c}}
  {\rho _{11}^r}& \cdots &{\rho _{1N}^r} \\
   \vdots & \ddots & \vdots  \\
  {\rho _{N1}^r}& \cdots &{\rho _{NN}^r}
\end{array}} \right],
\end{align}
\begin{align}
\rho _{nn'}^r = \exp \left( { - \varpi \Delta \phi _{nn'}^r}
\right),
\end{align}
\begin{align}
{{\bm{R}}^t} = \left[ {\begin{array}{*{20}{c}}
  {\rho _{11}^t}& \cdots &{\rho _{1M}^t} \\
   \vdots & \ddots & \vdots  \\
  {\rho _{M1}^t}& \cdots &{\rho _{MM}^t}
\end{array}} \right],
\end{align}
and
\begin{align}
\rho _{mm'}^t = \exp \left( { - \varpi \Delta \phi _{mm'}^t}
\right).
\end{align}
The symbol $\Delta\phi _{nn'}^r$ denotes the separation angle
between the $n$th and $n'$th transmitter-to-target paths,
$\Delta\phi_{mm'}^t$ denotes the separation angle between the $m$th
and $m'$th target-to-receiver paths, and $\varpi$ sets the
exponential decay in correlation with angle. From the model, it is
easy to see that larger $\varpi$ implies less dependency for fixed
$\Delta\phi _{nn'}^r$ and $\Delta\phi_{mm'}^t$. We consider $\varpi=
0.01$, $0.1$, and $\infty$ in the figures. Here $\varpi=\infty$
implies that the reflection coefficients are independent.


Figure \ref{MSEFANSHE} shows the comparison of RECRBOB and RMSE for
different $\varpi$ when all of the other parameters are the same as
in Figure \ref{5423}. We can see that the thresholds for the cases
with $\varpi=\infty$, $0.1$, and $0.01$ are 15 dB, 20 dB and 25 dB,
respectively. Thus, less dependency leads to a more favorable
threshold such that the RECRBOB is achievable at lower SCNR.  Above
threshold, all the curves are relatively close. The results imply that
the dependency of the reflection coefficients does not have
tremendous impact on the radar estimation performance, provided we
operate above threshold.  However, with
less dependency the radar can operate at lower SCNR while still
achieving an acceptable performance level.



\begin{figure}
\centering
 \includegraphics[width=3.0 in]{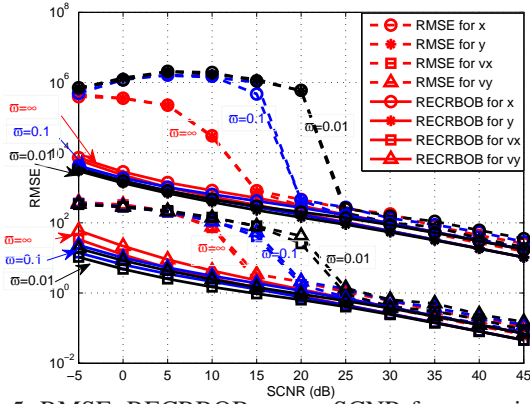}
\vspace{-0.4cm} \caption{RMSE, RECRBOB versus SCNR for a passive
MIMO radar with $M=2$ and $N=3$, spatially dependent reflection
coefficients, spatially independent noise, and orthogonal
signals.}\label{MSEFANSHE} \vspace{-0.1cm}
\end{figure}


\subsubsection{{Gaussian Spatially Dependent Noise}}

In this section, we consider the situation of spatially independent
reflection coefficients, spatially dependent noise, and orthogonal
signals. The elements of the noise covariance matrix $\bm{Q}$ are
generated with the following model
\begin{align}
{\bm{Q}} = \sigma _w^2{\bm{\tilde Q}} \otimes {{\bm{I}}_{K}},
\end{align}
where the $nn'$th element of $\bm{\tilde Q}$ is assumed to be
\begin{align}
{{\tilde Q}_{nn'}} = \exp \left\{ { - {d_{nn'}}\gamma } \right\}
\end{align}
where ${d_{nn'}}= \sqrt {{{\left( {x_n^r - x_{n'}^r} \right)}^2} +
{{\left( {y_n^r - y_{n'}^r} \right)}^2}}$, $\gamma$ sets the
exponential decay in correlation with distance, and ${{\bm{I}}_{K}}$
denotes a $K\times K$ identity matrix. From the model we can see
larger $\gamma$ results in less dependency for fixed $ d_{nn'} $. We
consider the situations of $\gamma  = 0.000005$, $0.00001$, and
$\infty$ and assume all the other parameters are the same as in
Figure \ref{5423}. Here $\gamma=\infty$ implies that the noise components are
independent.


Figure \ref{MSEFEIBAIZAO} shows the comparison of the RECRBOB and
RMSE for different values of $\gamma$.
{Case 1, Case 2 and Case 3
respectively represents  $\gamma=\infty$, $0.00001$, and
$0.000005$.} It is observed that the thresholds for cases with
$\gamma=\infty$, $0.00001$, and $ 0.000005$ are 15 dB, 10 dB, and 5
dB, respectively. Thus, more dependency leads to a more favorable
threshold such that the RECRBOB is achievable at lower SCNR. Above
the threshold, we see that $\gamma=0.000005$ has the smallest
RECRBOB while $\gamma=\infty$ has the largest RECRBOB, which means
larger dependency can lead to lower RECRBOB. In the cases considered
in Figure \ref{MSEFEIBAIZAO}, correlated noise leads to better
performance.


\begin{figure}
\centering
 \includegraphics[width=3.0 in]{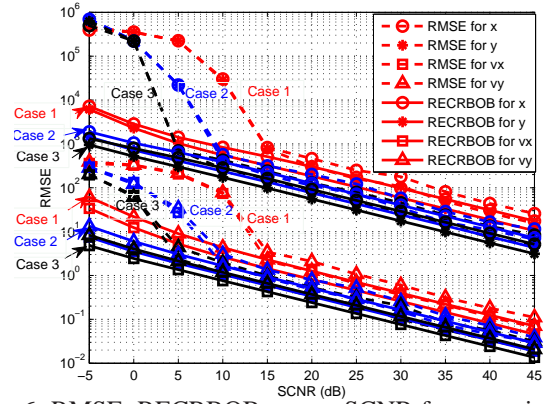}
\vspace{-0.4cm} \caption{RMSE, RECRBOB versus SCNR for a passive
MIMO radar with $M=2$ and $N=3$, spatially independent reflection
coefficients, dependent noise, and orthogonal
signals.}\label{MSEFEIBAIZAO} \vspace{-0.1cm}
\end{figure}



\subsubsection{{Inaccurate Signal Estimation}}

Now consider the case where the transmitted signals are not
estimated perfectly, possibly from the direct path receptions. Let
$n_{nm}(k),n=1,\cdots,N, m=1,\cdots,M, k=1,\cdots,K$ denote an
independent and identically distributed sequence of complex Gaussian
noise samples, each with zero mean and variance $0.1$ which models
the estimation error in the signal using
\begin{align}
{u_{nm}}(k) = \sqrt {\frac{{{E_m}{P_0}}}{{d_{tm}^2d_{rn}^2}}}
[{s_m}(k{T_s} - {\tau _{nm}},{\bm{\alpha} _m}) +
n{_{nm}}(k)]{e^{j2\pi {f_{nm}}k{T_s}}}\label{bbb}
\end{align}
Then (\ref{bbb}) is used to form $\bm{S}$ with the equations (9),
(10), (13) already given in the paper, and we call this mismatched
$\bm{S}$ $\bm{S_0}$. The undistorted $\bm{S}$ obtained this way, but
without additive noise, is called $\bm{S_1}$. This is exactly a case
where the model we employ in our estimation algorithm is mismatched
so the $\rm{RECRBOB_{mis}}$ results from Section~\ref{mm} become
applicable.  The resulting average $\rm{RECRBOB_{mis}}$ and
$\rm{RMSE_{mis}}$, after averaging over the noise using a Monte
Carlo simulation, are plotted in Figure \ref{MISCRB}.  From the
figure, we can see that the $\rm{RECRBOB_{mis}}$ provides an
informative lower bound\footnote{We have verified that the
unaveraged values of $\rm{RECRBOB_{mis}}$ also provide a lower bound
to the unaveraged values of $\rm{RMSE_{mis}}$.} on $\rm{RMSE_{mis}}$
in this case. Note that all the other details of the system analyzed
in Figure \ref{MISCRB}, except for this signal mismatch, are the
same as in Figure \ref{5423}.
\begin{figure}
\centering
 \includegraphics[width=3.0 in]{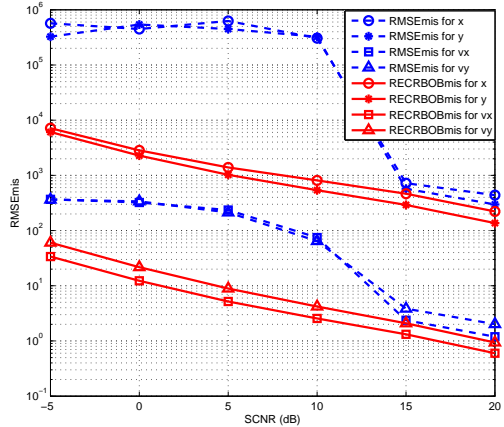}
\vspace{-0.4cm} \caption{RMSEmis, RECRBOBmis versus SCNR for a passive
MIMO radar with $M=2$ and $N=3$, spatially independent reflection
coefficients, independent noise, and orthogonal
signals.}\label{MISCRB} \vspace{-0.1cm}
\end{figure}

\section{Conclusions}\label{sec:Conlusion}

In this paper, we studied the performance of joint target position
and velocity estimation using a distributed radar network under more
general conditions than assumed in previous work. A received signal
model has been developed for active and passive radar with $M$
transmit and $N$ receive stations. The ML estimate and the exact CRB
expression are derived for possibly nonorthogonal signals, spatially
dependent Gaussian reflection coefficients, and spatially dependent
Gaussian clutter-plus-noise. {For cases in which some parameters
(for example the transmitted signal {from} direct path reception)
are not estimated correctly, we also derive the} {mismatched CRB.}
Numerical results are given to illustrate the use of the CRB and
mismatched CRB.  The numerical results show various cases with
signals of opportunity taken from a GSM wireless communication
system. It was shown that in the particular cases investigated, the
nonorthogonality of signal degraded the estimation performance both
in terms of RECRBOB and in terms of the threshold above which the
RMSE starts to become close to the RECRBOB in value and slope.
Decreasing the dependency between the different reflection
coefficients led to a more favorable threshold such that the radar
can operate at lower SCNR while still achieving an acceptable
performance level. Above threshold the dependency of the reflection
coefficients had little impact on the estimation performance,
provided a well performing estimation approach (nearly CRB
achieving) is employed. In some specific examples, it was also shown
that an increase in the dependency between the noise samples at
different antennas led to better estimation performance in terms of
the threshold and RECRBOB.

The work here can be generalized in several directions.
The CRB, a tight bound only for the high SCNR region
and being limited to unbiased estimators, is incapable of
characterizing the threshold value or accurately describing the
low-SCNR estimation performance of estimators. In this regard, we
need better analytical tools which can predict estimation
performance for low SCNR. The Ziv-Zakai bound is one promising
approach, which will be studied in our future work. The Ziv-Zakai bound
also allows prior information to be incorporated into the estimation.

\appendices


\section{Calculation of ${\bm{J}(\bm{\vartheta}|\bm{\alpha})}$} \label{app:computA}
According to {(23)} and (\ref{Jvarij}), we can obtain the FIM of the
vector $\bm{\vartheta}$ as
\begin{align}
\begin{array}{l}
{\bm{J}}\left( {\bm{\vartheta}|\bm{\alpha}}  \right)
 = \left[ \begin{array}{l}
{{\bm{J}}_{{\bm{\tau \tau }}}}{\kern 1pt} {\kern 1pt} {\kern 1pt} {\kern 1pt} {\kern 1pt} {\kern 1pt} {\kern 1pt} {\kern 1pt} {\kern 1pt} {\kern 1pt} {{\bm{J}}_{{\bm{\tau f}}}}{\kern 1pt} {\kern 1pt} {\kern 1pt} {\kern 1pt} {\kern 1pt} {\kern 1pt} {\kern 1pt} {\kern 1pt} {\kern 1pt} {\kern 1pt} {\kern 1pt} {\kern 1pt} {{\bm{J}}_{{\bm{\tau }}{{\bm{d}}_{\bm{t}}}}}{\kern 1pt} {\kern 1pt} {\kern 1pt} {\kern 1pt} {\kern 1pt} {\kern 1pt} {\kern 1pt} {\kern 1pt} {\kern 1pt} {\kern 1pt} {\kern 1pt} {{\bm{J}}_{{\bm{\tau }}{{\bm{d}}_{\bm{r}}}}}\\
{{\bm{J}}_{{\bm{f\tau }}}}{\kern 1pt} {\kern 1pt} {\kern 1pt} {\kern 1pt} {\kern 1pt} {\kern 1pt} {\kern 1pt} {\kern 1pt} {\kern 1pt} {\kern 1pt} {{\bm{J}}_{{\bm{ff}}}}{\kern 1pt} {\kern 1pt} {\kern 1pt} {\kern 1pt} {\kern 1pt} {\kern 1pt} {\kern 1pt} {\kern 1pt} {\kern 1pt} {\kern 1pt} {\kern 1pt} {\kern 1pt} {{\bm{J}}_{{\bm{f}}{{\bm{d}}_{\bm{t}}}}}{\kern 1pt} {\kern 1pt} {\kern 1pt} {\kern 1pt} {\kern 1pt} {\kern 1pt} {\kern 1pt} {\kern 1pt} {\kern 1pt} {\kern 1pt} {\kern 1pt} {{\bm{J}}_{{\bm{f}}{{\bm{d}}_{\bm{r}}}}}\\
{{\bm{J}}_{{{\bm{d}}_{\bm{t}}}{\bm{\tau }}}}{\kern 1pt} {\kern 1pt} {\kern 1pt} {\kern 1pt} {\kern 1pt} {\kern 1pt} {\kern 1pt} {\kern 1pt} {{\bm{J}}_{{{\bm{d}}_{\bm{t}}}{\bm{f}}}}{\kern 1pt} {\kern 1pt} {\kern 1pt} {\kern 1pt} {\kern 1pt} {\kern 1pt} {\kern 1pt} {{\bm{J}}_{{{\bm{d}}_{\bm{t}}}{{\bm{d}}_{\bm{t}}}}}{\kern 1pt} {\kern 1pt} {\kern 1pt} {\kern 1pt} {\kern 1pt} {\kern 1pt} {\kern 1pt} {\kern 1pt} {\kern 1pt} {{\bm{J}}_{{{\bm{d}}_{\bm{t}}}{{\bm{d}}_{\bm{r}}}}}\\
{{\bm{J}}_{{{\bm{d}}_{\bm{r}}}{\bm{\tau }}}}{\kern 1pt} {\kern 1pt}
{\kern 1pt} {\kern 1pt} {\kern 1pt} {\kern 1pt} {\kern 1pt} {\kern
1pt} {{\bm{J}}_{{{\bm{d}}_{\bm{r}}}{\bm{f}}}}{\kern 1pt} {\kern 1pt}
{\kern 1pt} {\kern 1pt} {\kern 1pt}
{{\bm{J}}_{{{\bm{d}}_{\bm{r}}}{{\bm{d}}_{\bm{t}}}}}{\kern 1pt}
{\kern 1pt} {\kern 1pt} {\kern 1pt} {\kern 1pt} {\kern 1pt} {\kern
1pt} {\kern 1pt} {\kern 1pt} {\kern 1pt} {\kern 1pt}
{{\bm{J}}_{{{\bm{d}}_{\bm{r}}}{{\bm{d}}_{\bm{r}}}}}
\end{array} \right]
\end{array}\label{abc},
\end{align}
where ${{\bm{J}}_{{\bm{\tau \tau }}}} = {\bm{J}}_{\bm{\tau
}}^H{{\bm{J}}_{\bm{\tau }}}$, ${{\bm{J}}_{{\bm{\tau f}}}}
={{\bm{J}}_{{\bm{f \tau}}}^H}= {\bm{J}}_{\bm{\tau
}}^H{{\bm{J}}_{\bm{f}}}$, ${{\bm{J}}_{{\bm{\tau
}}{{\bm{d}}_{\bm{t}}}}} = {{\bm{J}}_{{{\bm{d}}_{\bm{t}}}{\bm{\tau
}}}^H}={\bm{J}}_{\bm{\tau }}^H{{\bm{J}}_{{{\bm{d}}_{\bm{t}}}}}$,
${{\bm{J}}_{{\bm{\tau }}{{\bm{d}}_{\bm{r}}}}} =
{{\bm{J}}_{{{\bm{d}}_{\bm{r}}}{\bm{\tau }}}^H}={\bm{J}}_{\bm{\tau
}}^H{{\bm{J}}_{{{\bm{d}}_{\bm{r}}}}}$, ${{\bm{J}}_{{\bm{ff}}}} =
{\bm{J}}_{{\bm{f}}}^H{{\bm{J}}_{\bm{f}}}$,
${{\bm{J}}_{{\bm{f}}{{\bm{d}}_{\bm{t}}}}}
={{\bm{J}}_{{{\bm{d}}_{\bm{t}}}{\bm{f}}}^H}=
{\bm{J}}_{\bm{f}}^H{{\bm{J}}_{{{\bm{d}}_{\bm{t}}}}}$,
${{\bm{J}}_{{\bm{f}}{{\bm{d}}_{\bm{r}}}}}
={{\bm{J}}_{{{\bm{d}}_{\bm{r}}}{\bm{f}}}^H}=
{\bm{J}}_{\bm{f}}^H{{\bm{J}}_{_{{{\bm{d}}_{\bm{r}}}}}}$,
${{\bm{J}}_{{{\bm{d}}_{\bm{t}}}{{\bm{d}}_{\bm{t}}}}}{\kern 1pt}  =
{\bm{J}}_{{{\bm{d}}_{\bm{t}}}}^H{{\bm{J}}_{{{\bm{d}}_{\bm{t}}}}}$,
${{\bm{J}}_{{{\bm{d}}_{\bm{t}}}{{\bm{d}}_{\bm{r}}}}}
={{\bm{J}}_{{{\bm{d}}_{\bm{r}}}{{\bm{d}}_{\bm{t}}}}^H}=
{\bm{J}}_{{{\bm{d}}_{\bm{t}}}}^H{{\bm{J}}_{{{\bm{d}}_{\bm{r}}}}}$,
${{\bm{J}}_{{{\bm{d}}_{\bm{r}}}{{\bm{d}}_{\bm{r}}}}} =
{\bm{J}}_{{{\bm{d}}_{\bm{r}}}}^H{{\bm{J}}_{{{\bm{d}}_{\bm{r}}}}}$,
\begin{align}
{{\bm{J}}_{\bm{\tau }}} = \left( {{{\bm{C}}^{- \dag/2}} \otimes
{{\bm{C}}^{ - 1/2}}} \right)\frac{{\partial
{{\bm{C}}_{vec}}}}{{\partial {{\bm{\tau }}^\dag}}},
\end{align}
\begin{align}
{{\bm{J}}_{\bm{f}}} = \left( {{{\bm{C}}^{ - \dag/2}} \otimes
{{\bm{C}}^{ - 1/2}}} \right)\frac{{\partial
{{\bm{C}}_{vec}}}}{{\partial {{\bm{f}}^\dag}}},
\end{align}
\begin{align}
{{\bm{J}}_{{{\bm{d}}_{\bm{t}}}}} = \left( {{{\bm{C}}^{ - \dag/2}}
\otimes {{\bm{C}}^{ - 1/2}}} \right)\frac{{\partial
{{\bm{C}}_{vec}}}}{{\partial {{\bm{d}}_{\bm{t}}}^\dag}},
\end{align}
and
\begin{align}
{{\bm{J}}_{{{\bm{d}}_{\bm{r}}}}} = \left( {{{\bm{C}}^{ - \dag/2}}
\otimes {{\bm{C}}^{ - 1/2}}} \right)\frac{{\partial
{{\bm{C}}_{vec}}}}{{\partial {{\bm{d}}_{\bm{r}}}^\dag}}.
\end{align}
Then we reformulate the ${\bm{J}}\left( \bm{\vartheta} |\bm{\alpha}
\right)$ in a somewhat more explicit matrix form.

{First we derive $\bm{J}_{\bm{\tau\tau}}$, }let $\bm{s}_i$ and
$\bm{z}_i$ denote the $i$th column of $\bm S$ and $\bm{R}$,
respectively, such that ${\bm{S}} = \left[ {{{\bm{s}}_1}, \cdots
,{{\bm{s}}_{MN}}} \right]$ and ${\bm{R}} = \left[ {{{\bm{z}}_1},
\cdots ,{{\bm{z}}_{MN}}} \right]$. Note that $\bm{R}$ is a Hermitian
matrix, i.e., ${\bm{R}} = {\left[ {{\bm{z}}_1, \cdots
,{\bm{z}}_{MN}} \right]^H}$. Then, we have
\begin{align}
  \frac{{\partial {\bm{C}}}}{{\partial {\tau _{nm}}}} =& \frac{{\partial \left( {{\bm{SR}}{{\bm{S}}^H} + {\bm{Q}}} \right)}}{{\partial {\tau _{nm}}}} \notag\hfill \\
   =& \frac{{\partial {\bm{S}}}}{{\partial {\tau _{nm}}}}{\bm{R}}{{\bm{S}}^H} + {\bm{SR}}\frac{{\partial {{\bm{S}}^H}}}{{\partial {\tau _{nm}}}} \notag\hfill \\
   = &{\bm{s}}_i^{\tau}{\bm{z}}_i^H{{\bm{S}}^H} + {\bm{Sz}}_i^{}({\bm{s}}_i^{\tau})^H,
\end{align}
where $i=M(n-1)+m$ for $n=1,\cdots,N$ and $m=1,\cdots,M$,
\begin{align}
{\bm{s}}_i^\tau  =\frac{{\partial {{\bm{s}}_i}}}{{\partial {\tau
_{nm}}}} ={{\bm{e}}_n} \otimes {\left[ {\frac{{\partial
{u_{nm}}\left( 1 \right)}}{{\partial {\tau _{nm}}}}, \cdots
,\frac{{\partial {u_{nm}}\left( K \right)}}{{\partial {\tau
_{nm}}}}} \right]^\dag},
\end{align}
where $\bm{e}_n$ is {an} $N\times1$ column vector {with} zero
everywhere except for a 1 in the $n$th entry and
\begin{align}
\frac{{\partial {u_{nm}}(k)}}{{\partial {\tau _{nm}}}} =
{\sqrt{\frac{ {{E_m}{P_0}} }{{d^2_{tm}d^2_{rn}}}}}\frac{{\partial
{s_m}(k{T_s} - {\tau _{nm}},{{\bm{\alpha}}_m})}}{{\partial {\tau
_{nm}}}}{e^{j2\pi {f_{nm}}k{T_s}}}.
\end{align}

According to the following identity \cite{2012}
\begin{align}
\left( {{{\bm{X}}^\dag} \otimes {\bm{A}}} \right){\text{vec}}\left(
{\bm{B}} \right) = {\text{vec}}\left( {{\bm{ABX}}} \right),
\end{align}
we can obtain
\begin{align}
  {{\bm{J}}_{{\tau _{nm}}}} = &\left( {{{\bm{C}}^{ - \dag/2}} \otimes {{\bm{C}}^{ - 1/2}}} \right)\frac{{\partial {\text{vec}}\left( {\bm{C}} \right)}}{{\partial {\tau _{nm}}}}\notag \hfill \\
=& {\text{vec}}\left( {{{\bm{C}}^{ - 1/2}}\frac{{\partial {\bm{C}}}}{{\partial {\tau _{nm}}}}{{\bm{C}}^{ - 1/2}}} \right) \notag\hfill \\
= &{\text{ vec}}\left\{ {{{\bm{C}}^{ - 1/2}}\left( {{\bm{s}}_i^\tau {\bm{z}}_i^H{{\bm{S}}^H} + {\bm{Sz}}_i^{}{{({\bm{s}}_i^\tau )}^H}} \right){{\bm{C}}^{ - 1/2}}} \right\} \notag\hfill \\
= &{\text{vec}}\left( {{{\bm{V}}_i} + {\bm{V}}_i^H} \right)
\label{Jtaunm},
\end{align}
where ${{\bm{V}}_i} = {{\bm{C}}^{ - 1/2}}{\bm{s}}_i^\tau
{\bm{z}}_i^H{{\bm{S}}^H}{{\bm{C}}^{ - 1/2}}$. Using (\ref{Jtaunm})
and the following identity \cite{2012}
\begin{align}
{\text{vec}}{\left( {{{\bm{A}}^\dag}}
\right)^\dag}{\text{vec}}\left( {\bm{B}} \right) = Tr\left(
{{\bm{AB}}} \right) = Tr\left( {{\bm{BA}}} \right),
\end{align}
we can derive the $ij$th element of ${\bm{J}}^{\bm{\tau \tau}}$ as
follows {\small
\begin{align}
 & [{\bm{J}}_{\bm{\tau \tau}}]_{ij} ={\text{vec}}{\left\{ {{{\bm{V}}_i} + {\bm{V}}_i^H} \right\}^H}{\text{vec}}\left\{ {{{\bm{V}}_j} + {\bm{V}}_j^H} \right\} \notag\hfill \\
&= \text{Tr}\left\{ {\left( {{{\bm{V}}_i} + {\bm{V}}_i^H} \right)\left( {{{\bm{V}}_j} + {\bm{V}}_j^H} \right)} \right\}\notag \hfill \\
&= 2{\mathop{\Re}\nolimits} \{ Tr({V_i}{V_j} + V_i^H{V_j})\} \notag\hfill \\
&=2\Re \left\{ {{\bm{z}}_i^H{{\bm{S}}^H}{{\bm{C}}^{ -
1}}{\bm{s}}_j^\tau {\bm{z}}_j^H{{\bm{S}}^H}{{\bm{C}}^{ -
1}}{\bm{s}}_i^\tau + {{({\bm{s}}_i^\tau )}^H}{{\bm{C}}^{ -
1}}{\bm{s}}_j^\tau {\bm{z}}_j^H{{\bm{S}}^H}{{\bm{C}}^{ -
1}}{\bm{S}}{{\bm{z}}_i}} \right\}\label{Jtautauij},
\end{align}
} where $j=M(n'-1)+m'$ for $n'=1,\cdots,N$ and $m'=1,\cdots,M$.
Then, according to (\ref{Jtautauij}), we can reformulate
${{\bm{J}}_{{\bm{\tau \tau }}}}$ in the form of a matrix
\begin{align}
{{\bm{J}}_{{\bm{\tau \tau }}}} = 2{\Re} \left\{
{{\bm{Y}}{{\bm{S}}^{\bm{\tau}} } \odot {{\left(
{{\bm{Y}}{{\bm{S}}^{\bm{\tau}} }} \right)}^\dag} + {{\left(
{{{\bm{S}}^{\bm{\tau}} }} \right)}^H}{{\bm{C}}^{ -
1}}{{\bm{S}}^{\bm{\tau}} } \odot {{\left( {{\bm{YSR}}}
\right)}^\dag}} \right\}\label{Jtautau},
\end{align}
where ${{\bm{S}}^{\bm{\tau}} } = \left[ {{\bm{s}}_1^{\bm{\tau}} ,
\cdots ,{\bm{s}}_{MN}^{\bm{\tau}} } \right]$ and ${\bm{Y}} =
{\bm{R}}{{\bm{S}}^H}{{\bm{C}}^{ - 1}}$. Similarly, we can obtain
\begin{align}
{{\bm{J}}_{{\bm{\tau f}}}} = 2{\Re} \left\{ {{\bm{Y}}{{\bm{S}}^{\bm
f}} \odot {{\left( {{\bm{Y}}{{\bm{S}}^{\bm{\tau}} }} \right)}^\dag}
+ {{\left( {{{\bm{S}}^{\bm{\tau}} }} \right)}^H}{{\bm{C}}^{ -
1}}{{\bm{S}}^{\bm f}} \odot {{\left( {{\bm{YSR}}} \right)}^\dag}}
\right\},
\end{align}
and
\begin{align}
{{\bm{J}}_{{\bm{ff}}}} = 2{\Re} \left\{ {{\bm{Y}}{{\bm{S}}^{\bm f}}
\odot {{\left( {{\bm{Y}}{{\bm{S}}^{\bm f}}} \right)}^\dag} +
{{\left( {{{\bm{S}}^{\bm f}}} \right)}^H}{{\bm{C}}^{ -
1}}{{\bm{S}}^{\bm f}} \odot {{\left( {{\bm{YSR}}} \right)}^\dag}}
\right\}\label{Jff},
\end{align}
where ${{\bm{S}}^{\bm f }} = \left[ {{\bm{s}}_1^{\bm f} , \cdots
,{\bm{s}}_{MN}^{\bm{f}} } \right]$,
\begin{align}
{\bm{s}}_i^{\bm f}=\frac{{\partial {{\bm{s}}_i}}}{{\partial {f
_{nm}}}}={{\bm{e}}_n} \otimes {\left[ {\frac{{\partial
{u_{nm}}\left( 1 \right)}}{{\partial {f_{nm}}}}, \cdots
,\frac{{\partial {u_{nm}}\left( K \right)}}{{\partial {f_{nm}}}}}
\right]^\dag},
\end{align}
and
\begin{align}
\frac{{\partial {u_{nm}}\left( k \right)}}{{\partial {f_{nm}}}}
=j2\pi k{T_s}{u_{nm}}\left( k \right).
\end{align}

Next we derive ${{\bm{J}}_{{{\bm{d}}_{\bm{t}}}}}$,

\begin{align}
&\frac{{\partial {\bm{C}}}}{{\partial {d_{tm}}}} = \frac{{\partial ({\bm{SR}}{{\bm{S}}^H} + {\bm{Q}})}}{{\partial {d_{tm}}}}\notag\\
{\kern 1pt}  &= \frac{{\partial {\bm{S}}}}{{\partial {d_{tm}}}}{\bm{R}}{{\bm{S}}^H} + {\bm{SR}}\frac{{\partial {{\bm{S}}^H}}}{{\partial {d_{tm}}}}\notag\\
 &= ({\bm{s}}_m^t{\bm{z}}_m^H + {\bm{s}}_{m + M}^t{\bm{z}}_{m + M}^H +  \cdots  + {\bm{s}}_{m + (N - 1)M}^t{\bm{z}}_{m + (N - 1)M}^H){{\bm{S}}^H}\notag\\
 &+ {\bm{S}}({{\bm{z}}_m}{({\bm{s}}_m^t)^H} + {{\bm{z}}_{m + M}}{({\bm{s}}_{m + M}^t)^H} +  \cdots  + {{\bm{z}}_{m + (N - 1)M}}{({\bm{s}}_{m + (N -
 1)M}^t)^H}),
\end{align}
where
\begin{align}
&{\bm{s}}_{m + (n - 1)M}^t = \frac{{\partial {{\bm{s}}_{m + (n - 1)M}}}}{{\partial {d_{tm}}}}\notag\\
 &= {{\bm{e}}_n} \otimes {[\frac{{\partial {u_{nm}}(1)}}{{\partial {d_{tm}}}},\frac{{\partial {u_{nm}}(2)}}{{\partial {d_{tm}}}} \cdots ,\frac{{\partial {u_{nm}}(K)}}{{\partial {d_{tm}}}}]^\dag},{\kern 1pt} {\kern 1pt} n = 1, \cdots ,N
\end{align}
and
\begin{align}
\frac{{\partial {u_{nm}}(k)}}{{\partial {d_{tm}}}} =  - \frac{{\sqrt
{{E_m}{P_0}} }}{{d_{tm}^2d_{rn}}}s(k{T_s} - {\tau
_{nm}},{\bm{\alpha} _m}){e^{j2\pi {f_{nm}}k{T_s}}}.
\end{align}

It can be derived that
\begin{align}
{{\bm{J}}_{{d_{tm}}}} = &({{\bm{C}}^{ - \dag/2}} \otimes {{\bm{C}}^{ - 1/2}})\frac{{\partial {{\bm{C}}_{vec}}}}{{\partial {d_{tm}}}}\notag\\
=&vec({{\bm{C}}^{ - 1/2}}\frac{{\partial {\bm{C}}}}{{\partial {d_{tm}}}}{{\bm{C}}^{ - 1/2}})\notag\\
=&vec\{ {{\bm{C}}^{ - 1/2}}(({\bm{s}}_m^t{\bm{z}}_m^H + {\bm{s}}_{m + M}^t{\bm{z}}_{m + M}^H +  \cdots \notag\\
 &+ {\bm{s}}_{m + (N - 1)M}^t{\bm{z}}_{m + (N - 1)M}^H){{\bm{S}}^H}  \notag\\
&+{\bm{S}}({{\bm{z}}_m}{({\bm{s}}_m^t)^H} + {{\bm{z}}_{m + M}}{({\bm{s}}_{m + M}^t)^H} +  \cdots \notag\\
 &+ {{\bm{z}}_{m + (N - 1)M}}{({\bm{s}}_{m + (N - 1)M}^t)^H})){{\bm{C}}^{ - 1/2}}\} \notag\\
 =&vec({\ell _m} + \ell _m^H),
\end{align}
where ${{\bm{\ell}} _m} = {{\bm{C}}^{-1/2}}({\bm{s}}_m^t{\bm{z}}_m^H
+ {\bm{s}}_{m + M}^t{\bm{z}}_{m + M}^H +  \cdots  + {\bm{s}}_{m + (N
- 1)M}^t{\bm{z}}_{m + (N - 1)M}^H){{\bm{S}}^H}{{\bm{C}}^{ - 1/2}}$.
Then, we obtain
\begin{align}
&{[{{\bm{J}}_{{{\bm{d}}_{\bm{t}}}{{\bm{d}}_{\bm{t}}}}}]_{mm'}} = vec{({\ell _m} + \ell _m^H)^H}vec({\ell _{m'}} + \ell _{m'}^H)\notag\\
 &= Tr\{ ({\ell _m} + \ell _m^H)({\ell _{m'}} + \ell _{m'}^H)\} \notag\\
 &= 2\Re \{ Tr({\ell _m}{\ell _{m'}} + \ell _m^H{\ell _{m'}})\} \notag\\
 &= 2\Re \{ \sum\limits_{n = 1}^N {\sum\limits_{n' = 1}^N {({{({{\bm{z}}_{m + (n - 1)M}})}^H}{{\bm{S}}^H}{{\bm{C}}^{ - 1}}{\bm{s}}_{m' + (n' - 1)M}^t{{({{\bm{z}}_{m' + (n' - 1)M}})}^H}{{\bm{S}}^H}} } \notag\\
&{{\bm{C}}^{ - 1}}{\bm{s}}_{m + (n - 1)M}^t + {({\bm{s}}_{m + (n - 1)M}^t)^H}{{\bm{C}}^{ - 1}}{\bm{s}}_{m' + (n' - 1)M}^t{({{\bm{z}}_{m' + (n' - 1)M}})^H}{{\bm{S}}^H}\notag\\
&{{\bm{C}}^{ - 1}}{\bm{S}}{{\bm{z}}_{m + (n - 1)M}})\}.
\end{align}

Reformulate ${{\bm{J}}_{{{\bm{d}}_{\bm{t}}}{{\bm{d}}_{\bm{t}}}}}$ in
the form of a matrix
\begin{align}
&{{\bm{J}}_{{{\bm{d}}_{\bm{t}}}{{\bm{d}}_{\bm{t}}}}} = 2\Re \{ \sum\limits_{n = 1}^N {\sum\limits_{n' = 1}^N {({\aleph _n}} } {{\bm{S}}^H}{{\bm{C}}^{ - 1}}{\Im _{n'}} \odot {({\aleph _{n'}}{{\bm{S}}^H}{{\bm{C}}^{ - 1}}{\Im _n})^\dag }\notag\\
& + {({\Im _n})^H}{{\bm{C}}^{ - 1}}{\Im _{n'}} \odot {({\aleph
_{n'}}{{\bm{S}}^H}{{\bm{C}}^{ - 1}}{\bm{S}}{({\aleph _n})^H})^\dag
})\},
\end{align}
where ${\aleph _n} = {({{\bm{z}}_{1 + (n - 1)M}}, \cdots
,{{\bm{z}}_{M + (n - 1)M}})^H}$, ${\Im _n} = ({\bm{s}}_{1 + (n -
1)M}^t, \cdots ,{\bm{s}}_{M + (n - 1)M}^t)$. Similarly, we can
derive
\begin{align}
&{{\bm{J}}_{{{\bm{d}}_{\bm{t}}}{\bm{\tau }}}} = 2\Re \{ \sum\limits_{n = 1}^N {({\aleph _n}} {{\bm{S}}^H}{{\bm{C}}^{ - 1}}{{\bm{S}}^{\bm{\tau }}} \odot {({\bm{R}}{{\bm{S}}^H}{{\bm{C}}^{ - 1}}{\Im _n})^\dag }\notag\\
& + {({\Im _n})^H}{{\bm{C}}^{ - 1}}{{\bm{S}}^{\bm{\tau }}} \odot
{({\bm{R}}{{\bm{S}}^H}{{\bm{C}}^{ - 1}}{\bm{S}}{({\aleph
_n})^H})^\dag })\},
\end{align}
\begin{align}
&{{\bm{J}}_{{{\bm{d}}_{\bm{t}}}{\bm{f}}}} = 2\Re \{ \sum\limits_{n = 1}^N {({\aleph _n}} {{\bm{S}}^H}{{\bm{C}}^{ - 1}}{{\bm{S}}^{\bm{f}}} \odot {({\bm{R}}{{\bm{S}}^H}{{\bm{C}}^{ - 1}}{\Im _n})^\dag }\notag\\
& + {({\Im _n})^H}{{\bm{C}}^{ - 1}}{{\bm{S}}^{\bm{f}}} \odot
{({\bm{R}}{{\bm{S}}^H}{{\bm{C}}^{ - 1}}{\bm{S}}{({\aleph
_n})^H})^\dag })\}.
\end{align}

{To derive $d_{rn}$, we employ }
\begin{align}
&\frac{{\partial C}}{{\partial {d_{rn}}}} = \frac{{\partial ({\bm{SR}}{{\bm{S}}^H} + {\bm{Q}})}}{{\partial {d_{rn}}}}\notag\\
 &= \frac{{\partial {\bm{S}}}}{{\partial {d_{rn}}}}{\bm{R}}{{\bm{S}}^H} + {\bm{SR}}\frac{{\partial {{\bm{S}}^H}}}{{\partial {d_{rn}}}}\notag\\
 &= ({\bm{s}}_{1 + (n - 1)M}^r{\bm{z}}_{1 + (n - 1)M}^H + {\bm{s}}_{2 + (n - 1)M}^r{\bm{z}}_{1 + (n - 1)M}^H +  \cdots \notag\\
 &+ {\bm{s}}_{M + (n - 1)M}^r{\bm{z}}_{M + (n - 1)M}^H){{\bm{S}}^H}\notag\\
&{\kern 1pt}  + {\bm{S}}({{\bm{z}}_{1 + (n - 1)M}}{(s_{1 + (n - 1)M}^r)^H} + {{\bm{z}}_{2 + (n - 1)M}}{({\bm{s}}_{2 + (n - 1)M}^t)^H} +  \cdots  + {{\bm{z}}_{M + (n - 1)M}}\notag\\
&{({\bm{s}}_{M + (n - 1)M}^t)^H}),
\end{align}
where
\begin{align}
&{\bm{s}}_{m + (n - 1)M}^r = \frac{{\partial {{\bm{s}}_{m + (n - 1)M}}}}{{\partial {d_{rn}}}}\notag\\
 &= {{\bm{e}}_n} \otimes {[\frac{{\partial {u_{nm}}(1)}}{{\partial {d_{rn}}}},\frac{{\partial {u_{nm}}(2)}}{{\partial {d_{rn}}}} \cdots ,\frac{{\partial {u_{nm}}(K)}}{{\partial {d_{rn}}}}]^\dag},{\kern 1pt} {\kern 1pt} {\kern 1pt} m = 1, \cdots
 ,M,
\end{align}
and
\begin{align}
\frac{{\partial {u_{nm}}(k)}}{{\partial {d_{rn}}}} =  - \frac{{\sqrt
{{E_m}{P_0}} }}{{d_{tm}d_{rn}^2}}s(k{T_s} - {\tau
_{nm}},{{\bm{\alpha }}_m}){e^{j2\pi {f_{nm}}k{T_s}}},
\end{align}

We can {then} derive
\begin{align}
&{{\bm{J}}_{{d_{rn}}}} = ({{\bm{C}}^{ - \dag/2}} \otimes {{\bm{C}}^{ - 1/2}})\frac{{\partial {{\bm{C}}_{vec}}}}{{\partial {d_{rn}}}}\notag\\
&{\kern 1pt}  = vec({{\bm{C}}^{ - 1/2}}\frac{{\partial {\bm{C}}}}{{\partial {d_{rn}}}}{{\bm{C}}^{ - 1/2}})\notag\\
 &= vec\{ {{\bm{C}}^{ - 1/2}}(({\bm{s}}_{1 + (n - 1)M}^r{\bm{z}}_{1 + (n - 1)M}^H + {\bm{s}}_{2 + (n - 1)M}^r{\bm{z}}_{2 + (n - 1)M}^H +  \cdots \notag\\
 &+ {\bm{s}}_{M + (n - 1)M}^r{\bm{z}}_{M + (n - 1)M}^H){{\bm{S}}^H}\notag\\
 &+ {\bm{S}}({({\bm{s}}_{1 + (n - 1)M}^r)^H}{{\bm{z}}_{1 + (n - 1)M}} + {({\bm{s}}_{2 + (n - 1)M}^t)^H}{{\bm{z}}_{2 + (n - 1)M}} +  \cdots \notag\\
 &+ {({\bm{s}}_{M + (n - 1)M}^t)^H}{{\bm{z}}_{M + (n - 1)M}})){{\bm{C}}^{ - 1/2}}\} \notag\\
&= vec({w_n} + w_n^H),
\end{align}
where ${w_n} = {{\bm{C}}^{ - 1/2}}({\bm{s}}_{1 + (n -
1)M}^r{\bm{z}}_{1 + (n - 1)M}^H + {\bm{s}}_{2 + (n -
1)M}^r{\bm{z}}_{2 + (n - 1)M}^H +  \cdots  + {\bm{s}}_{M + (n -
1)M}^r{\bm{z}}_{M + (n - 1)M}^H){{\bm{S}}^H}{{\bm{C}}^{ - 1/2}}$.

Then, we can obtain
\begin{align}
&{[{{\bm{J}}_{{{\bm{d}}_{\bm{r}}}{{\bm{d}}_{\bm{r}}}}}]_{nn'}} = vec{({w_n} + w_n^H)^H}vec({w_{n'}} + \ell _{n'}^H)\notag\\
 &= Tr\{ ({w_n} + w_n^H)({w_{n'}} + w_{n'}^H)\} \notag\\
 &= 2\Re \{ Tr({w_n}{w_{n'}} + w_n^H{w_{n'}})\} \notag\\
 &= 2\Re \{ \sum\limits_{m = 1}^M {\sum\limits_{m' = 1}^M {({{({{\bm{z}}_{m + (n - 1)M}})}^H}{{\bm{S}}^H}{{\bm{C}}^{ - 1}}{\bm{s}}_{m' + (n' - 1)M}^t{{({{\bm{z}}_{m' + (n' - 1)M}})}^H}{{\bm{S}}^H}} } \notag\\
&{{\bm{C}}^{ - 1}}{\bm{s}}_{m + (n - 1)M}^t + {({\bm{s}}_{m + (n - 1)M}^t)^H}{{\bm{C}}^{ - 1}}{\bm{s}}_{m' + (n' - 1)M}^t{({{\bm{z}}_{m' + (n' - 1)M}})^H}{{\bm{S}}^H}\notag\\
&{{\bm{C}}^{ - 1}}{\bm{S}}{{\bm{z}}_{m + (n - 1)M}})\}\label{aaaaa}.
\end{align}
{The result of (\ref{aaaaa}) can be reformulated as  }
\begin{align}
&{{\bm{J}}_{{{\bm{d}}_{\bm{r}}}{{\bm{d}}_{\bm{r}}}}} = 2\Re \{ \sum\limits_{m = 1}^M {\sum\limits_{m' = 1}^M {({\mathchar'26\mkern-10mu\lambda _m}} } {{\bm{S}}^H}{{\bm{C}}^{ - 1}}{\wp _{m'}} \odot {({\mathchar'26\mkern-10mu\lambda _{m'}}{{\bm{S}}^H}{{\bm{C}}^{ - 1}}{\wp _m})^\dag }\notag\\
 &+ {({\wp _m})^H}{{\bm{C}}^{ - 1}}{\wp _{m'}} \odot {({\mathchar'26\mkern-10mu\lambda _{m'}}{{\bm{S}}^H}{{\bm{C}}^{ - 1}}{\bm{S}}{({\mathchar'26\mkern-10mu\lambda _m})^H})^\dag
 })\},
\end{align}
where ${\mathchar'26\mkern-10mu\lambda _m} = {({{\bm{z}}_m}, \cdots
,{{\bm{z}}_{m + (N - 1)M}})^H}$, ${\wp _m} = (s_{m}^r, \cdots ,s_{m
+ (N - 1)M}^r)$. Similarly, we can obtain
\begin{align}
&{{\bm{J}}_{{{\bm{d}}_{\bm{r}}}{\bm{\tau }}}} = 2\Re \{ \sum\limits_{m = 1}^M {({\mathchar'26\mkern-10mu\lambda _m}} {{\bm{S}}^H}{{\bm{C}}^{ - 1}}{{\bm{S}}^{\bm{\tau }}} \odot {({\bm{R}}{{\bm{S}}^H}{{\bm{C}}^{ - 1}}{\wp _m})^\dag }\notag\\
 &+ {({\wp _m})^H}{{\bm{C}}^{ - 1}}{{\bm{S}}^\tau } \odot {({\bm{R}}{{\bm{S}}^H}{{\bm{C}}^{ - 1}}{\bm{S}}{({\mathchar'26\mkern-10mu\lambda _m})^H})^\dag })\}
\end{align}
\begin{align}
&{{\bm{J}}_{{{\bm{d}}_{\bm{r}}}{\bm{f}}}} = 2\Re \{ \sum\limits_{m = 1}^M {({\mathchar'26\mkern-10mu\lambda _m}} {{\bm{S}}^H}{{\bm{C}}^{ - 1}}{{\bm{S}}^{\bm{f}}} \odot {({\bm{R}}{{\bm{S}}^H}{{\bm{C}}^{ - 1}}{\wp _m})^\dag }\notag\\
 &+ {({\wp _m})^H}{{\bm{C}}^{ - 1}}{{\bm{S}}^f} \odot {({\bm{R}}{{\bm{S}}^H}{{\bm{C}}^{ - 1}}{\bm{S}}{({\mathchar'26\mkern-10mu\lambda _m})^H})^\dag
 })\},
\end{align}
\begin{align}
&{{\bm{J}}_{{{\bm{d}}_{\bm{r}}}{{\bm{d}}_{\bm{t}}}}} = 2\Re \{ \sum\limits_{m = 1}^M {\sum\limits_{n = 1}^N {({\mathchar'26\mkern-10mu\lambda _m}} } {{\bm{S}}^H}{{\bm{C}}^{ - 1}}{\Im _n} \odot {({\aleph _n}{{\bm{S}}^H}{{\bm{C}}^{ - 1}}{\wp _m})^\dag }\notag\\
 &+ {({\wp _m})^H}{{\bm{C}}^{ - 1}}{\Im _n} \odot {({\aleph _n}{{\bm{S}}^H}{{\bm{C}}^{ - 1}}{\bm{S}}{({\mathchar'26\mkern-10mu\lambda _m})^H})^\dag })\}
\end{align}

\bibliographystyle{unsrt}

\end{document}